\documentclass{amsart}

\usepackage{amssymb}
\usepackage[dvips]{epsfig}
\usepackage{graphicx}
\usepackage{color}

 \usepackage[latin1]{inputenc}
 \usepackage[T1]{fontenc}
 \usepackage[normalem]{ulem}
\usepackage{verbatim}
 \usepackage{graphicx}
\usepackage[latin1]{inputenc}
\usepackage{amsmath}
\usepackage{amssymb}
\usepackage{amsfonts}
\usepackage{amsthm}
\usepackage{bbm}
 \usepackage{comment}

\usepackage{color}

\newcommand{\dsp}{\displaystyle}
\newcommand{\dt}{\partial_t}

\newcommand{\R}{{\mathbb R}}
\newcommand{\T}{{\mathbb T}}

\newcommand{\eps}{\varepsilon}

\newcommand{\D}{\vert D\vert}

\theoremstyle{remark}

\newtheorem{remark}{Remark}

\theoremstyle{definition}

\theoremstyle{theorem}

\newtheorem{lemma}{Lemma}
\newtheorem{theorem}{Theorem}
\newtheorem{corollary}{Corollary}

\newtheorem*{merci}{Acknowledgements}

\numberwithin{equation}{section}

\begin{document}

\title[BO-ILW]{ Benjamin-Ono and Intermediate Long Wave equations : modeling, IST and PDE}

\author[J.-C. Saut]{Jean-Claude Saut}
\address{Laboratoire de Math\' ematiques, UMR 8628\\
Universit\' e Paris-Saclay, Paris-Sud  et CNRS\\ 91405 Orsay, France}
\email{jean-claude.saut@u-psud.fr}

\date{November 20, 2018}
\maketitle

\begin{abstract}
This survey article is focused on two asymptotic models for internal waves, the Benjamin-Ono (BO) and Intermediate Long Wave  (ILW) equations that are integrable by inverse scattering techniques (IST). After recalling briefly their (rigorous) derivations we will review old and recent results on the Cauchy problem, comparing those obtained by IST and PDE techniques and also results more connected to the physical origin of the equations. We will consider mainly the Cauchy problem on the whole real line with only a few comments on the periodic case. We will also briefly discuss some close relevant problems in particular the higher order extensions and the two-dimensional (KP like) versions of the BO and ILW equations.

\end{abstract}

\section{Introduction}

In order to illustrate the links and interactions between PDE and Inverse Scattering methods, we have chosen to focus on two one-dimensional examples that have a  physical relevance (in the context of internal waves)  and that lead to yet unsolved interesting issues.

The  Benjamin-Ono  (BO) and Intermediate Long Wave (ILW)  equations are two classical examples of completely integrable one-dimensional equations, maybe not so well-known as the Korteweg de Vries or the cubic nonlinear Schr\"{o}dinger equations though. A striking fact is that a complete rigorous resolution of the Cauchy problem by IST techniques is still incomplete for the BO equations for arbitrary large initial data while it is open in the ILW case, even for small data. On the other hand, both those problems can be solved by "PDE" techniques, for arbitrary initial data in relatively big spaces but no general result on the long time behavior  of "large" solutions is known with the notable exception of stability issues of solitons and multisolitons.

Those two scalar equations are one-dimensional, one-way propagation asymptotic models for internal waves in an appropriate regime. We have thus three different viewpoints on BO and ILW equations and this article aims to review them and emphasize their possible links. Since we do not want to ignore the modeling aspects, we first recall the derivation of the equations in the context of internal waves in a two-layer system. The modeling of internal waves displays a variety of fascinating scientific problems. We refer for instance  to the survey article \cite{HM} for the physical modeling aspects and to \cite{BLS, CGK, Sa, OI} for the rigorous derivation of asymptotic models.

The  Benjamin-Ono equation was first formally derived by Brooke Benjamin, \cite{Ben} (and independently in  \cite{DaAc} where one can find also numerical simulations and experimental comparisons), and later by Ono \cite{O}, to describe the propagation of long weakly nonlinear internal waves in a stratified fluid such that two layers of different densities are joined by a thin region where the density varies continuously (pycnocline), the lower layer  being infinite. \footnote{One can find interesting comparisons with experiments in \cite{Max}.}
Benjamin also wrote down the explicit algebraically decaying solitary wave solution and also the periodic traveling wave.  
 

The Intermediate Long Wave equation was introduced  by Kubota {\it etal}  \cite{KKD} to describe the propagation of a long weakly nonlinear internal wave in a stratified medium of finite total depth.  A formal derivation   is also  given in Joseph \cite{J}  who used the dispersion relation derived in \cite{Phi} in the context of the Whitham non local equation \cite{Whi}. Joseph derived furthermore the solitary wave solution. 

We also refer to  \cite{MaRe, KB, CMH, LHA, Roma} for the relevance of the ILW equation in various oceanic or atmospheric contexts.

The ILW equation reduces formally to the BO equation when the depth of the lower layer tends to infinity. 

A rigorous derivation (in the sense of consistency) is given in \cite{BLS, CGK} using a two-layer system, that is a system of two layers of fluids of different densities, the density of the total fluid being discontinuous though (see below). Interesting comparisons with experiments can be found {\it eg} in \cite{KB}. \footnote{Recall however that BO and ILW equations are weakly nonlinear models and they do not fit well with the modeling of higher amplitude waves, see {\it eg} the experiments in \cite{SJ}.}Bidirectional versions are derived in \cite{ChCa, ChCa2, CGK, BLS}, see below for a quick description of a  rigorous derivation in the sense of consistency.

Both equations  belong to the general class of equations of the type (see \cite{LPS})

\begin{equation}\label{dBurg}
u_t+uu_x-\mathcal L u_x=0,
\end{equation}
where $\mathcal L$ is defined after Fourier transform by $\widehat{\mathcal L f}(\xi)=p(\xi)\hat f(\xi)$ where $p$ is  a real symbol, 

with $p(\xi)=p_\delta(\xi)=\xi\coth(\delta \xi)-\frac{1}{\delta}, \delta >0$ for the ILW equation and $p(\xi)=|\xi|$ for the BO equation.

They can alternatively be written respectively,

\begin{equation}\label{BO}
u_t+uu_x-Hu_{xx}=0,
\end{equation}

where H is the Hilbert transform,  that is the convolution with PV($\frac{1}{x}$).

and

\begin{equation}\label{ILW}
u_t+uu_x+\frac{1}{\delta}u_x+\mathcal T(u_{xx})=0,
\end{equation}

where 

$$\mathcal T=PV\int_{-\infty}^\infty \coth(\frac{x-y}{\delta})u(y)dy.$$

The BO equation has the following scaling and translation invariance :

if u is a solution de BO, so is v defined by 

$$v(x,t)=cu(c(x-x_0), c^2t)), \; \forall c>0, x_0\in \R.$$

As aforementioned, the ILW equation reduces formally to the BO equation when $\delta \to \infty$ and it was actually proven in \cite{ABFS} that the solution $u_\delta$ of \eqref{ILW} with initial data $u_0$ converges as $\delta \to +\infty$ to the solution of the Benjamin-Ono equation \eqref{BO} with the same initial data in suitable Sobolev spaces.

Furthermore, if $u_\delta$ is a solution of \eqref{ILW} and setting 

$$v_\delta(x,t)=\frac{3}{\delta}u_\delta(x,\frac{3}{\delta}t),$$

$v_\delta$ tends as $\delta \to 0$ to the solution $u$
of the KdV equation

\begin{equation}\label{KdVnew}
u_t+uu_x+u_{xxx}=0
\end{equation}

Both the BO and ILW equations conserve formally the $L^2$ norm of the initial data. Moreover they have an Hamiltonian structure

$$u_t+\partial_x\mathcal H u=0,$$

where  
 $$\mathcal H u=\frac{1}{2}\int_\R(\frac{u^3}{3}-(\vert D\vert^{1/2}u)^2 )dx$$

 for the BO equation and 
 
 $$\mathcal H u=\frac{1}{2}\int_\R(\frac{u^3}{3}-\frac{1}{2}\vert T_\delta^{1/2}u\vert^2),$$

 where

 $$\widehat{T_\delta f}(\xi)=p_\delta(\xi)=\xi\coth(\delta \xi)-\frac{1}{\delta},$$
 
 for the ILW equation.

\vspace{0.3cm}
 
 We now recall now the long process leading to ILW and BO equations following \cite{BLS} that will lead to a rigorous justification of the equations in the sense of consistency.

The physical context is that of the two-layer system of two inviscid incompressible fluids of different densities $\rho_1<\rho_2,$ see picture below.

\begin{center}
\includegraphics[width=.85\textwidth]{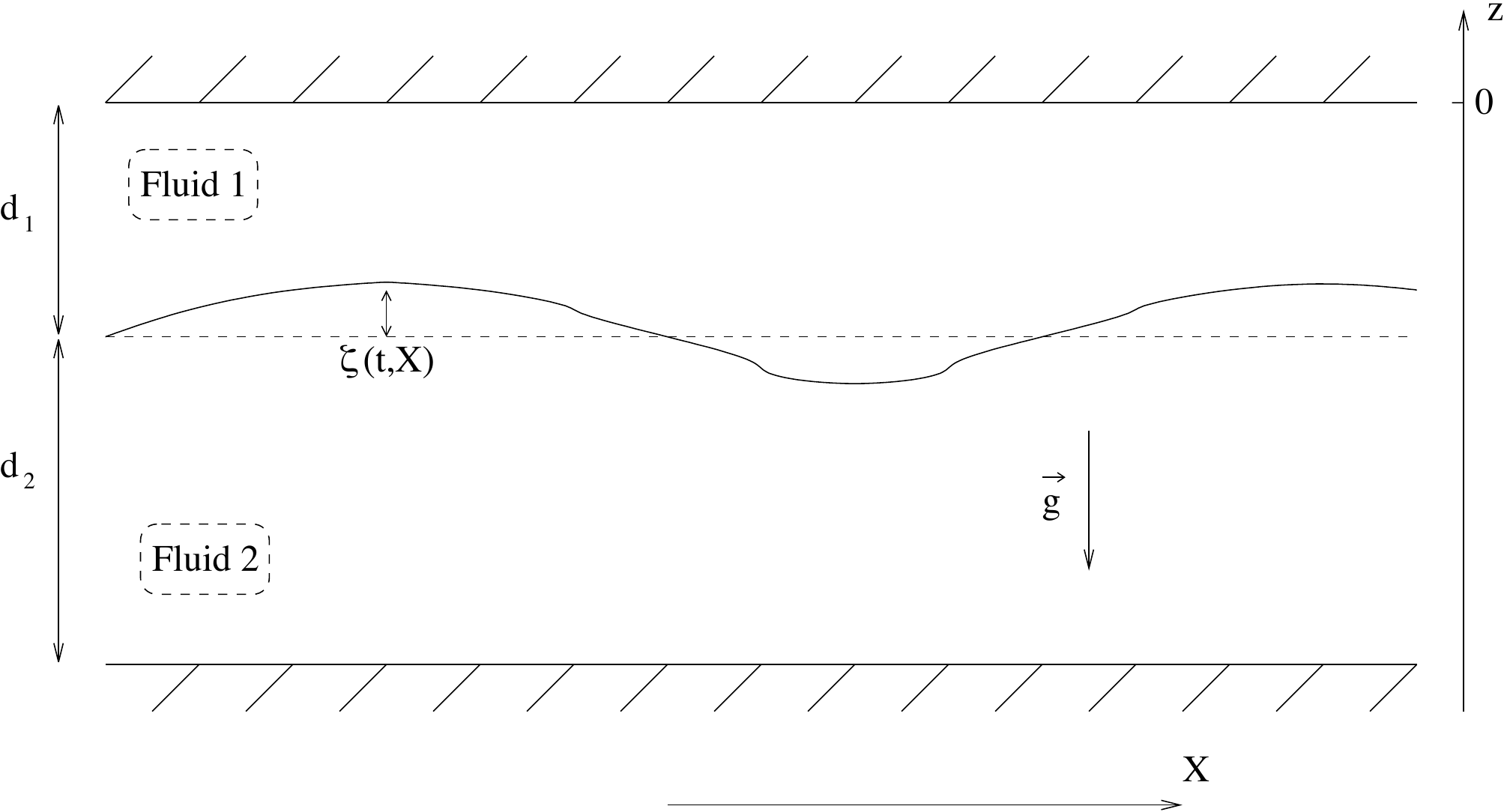}

\end{center}

\vspace{1cm}
We refer to \cite{BLS} for the derivation of the complete two-layer system (extending the classical water wave system, see \cite{La}) describing the evolution of the interface $\zeta$ and of a suitable velocity variable and to \cite{Lannes2} for  a deep analysis of this system.

The asymptotic models are derived from the two-layer system \footnote{Both the BO and ILW equations were first derive  in the context of a continuously stratified fluid to model long internal waves on a pycnocline  (boundary separating two liquid layers of different densities), see \cite{Ben, KKD, DaAc}). The rigorous justification in this context is studied in \cite{DLS}.}\cite{BLS} (see also \cite{Sa}) after introducing scaling parameters, namely

\begin{itemize}
\item   ${\bf a} =$  typical amplitude of the deformation of the
interface,   $ \lambda =$ typical wavelength.
 
\item Dimensionless independent variables
$$
    \widetilde{X}:=\frac{X}{\lambda},\quad
    \widetilde{z}:=\frac{z}{d_1},\quad
    \widetilde{t}:=\frac{t}{\lambda/\sqrt{gd_1},\quad},
$$
are introduced.  Likewise, we define the dimensionless unknowns
$$
    \widetilde{\zeta}:=\frac{\zeta}{a},\quad
    \widetilde{\psi}_1:=\frac{\psi_1}{a\lambda\sqrt{g/d_1}},
$$
as well as the dimensionless parameters
$$
    \gamma:=\frac{\rho_1}{\rho_2},\quad
    d:=\frac{d_1}{d_2},\quad
    \epsilon:=\frac{a}{d_1},\quad
    \mu:=\frac{d_1^2}{\lambda^2}; \footnote{Note that $d\sim \frac{1}{\delta},  $  where $\delta$ is as in the  the previous notation.}
$$
\end{itemize}
Though they are redundant, it is also notationally convenient to
introduce two other parameters $\eps_2$ and $\mu_2$ defined as
$$
    \eps_2=\frac{a}{d_2}=\eps d,\qquad
    \mu_2=\frac{d_2^2}{\lambda^2}=\frac{\mu}{d^2}.
$$

The range of validity of the various regimes is summarized in the
following table.

\vspace{0.5cm}
{\footnotesize
\begin{center}
\begin{tabular}{|*{3}{l|}}
\hline
     & $\epsilon=O(1)$ & $\epsilon\ll1$   \\
    \hline
    $\mu=O(1)$ & Full equations & $d\sim 1$: FD/FD eq'ns\\
    \hline
    $\mu\ll1$ & $d\sim 1$: SW/SW eq'ns & $\mu\sim\epsilon$ and $d^2\sim\epsilon$: B/FD eq'ns\\
    &  $d^2\sim \mu\sim\epsilon_2^2$: SW/FD eq'ns & $\mu\sim\epsilon$ and $d\sim1$: B/B eq'ns\\
    &  &$d^2\sim \mu\sim\epsilon^2$: ILW eq'ns \\
    &  &$d=0$ and $\mu\sim\epsilon^2$: BO eq'ns  \\
    \hline
    \end{tabular}
\end{center}
 }

\vspace{0.5cm}

The ILW regime is thus obtained when 
    $\mu\sim \eps^2\ll 1$ and $\mu_2\sim 1$
    (and thus $d^2\sim \mu\sim\eps_2$); in this case, one gets the following expansion of the nonlocal {\it interface operator} (see \cite{BLS} for details), where $\vert D\vert =\sqrt{-\Delta}$ :
    \begin{equation}\label{asympILW}
    {\bf H}^{\mu,d}[\eps\zeta]\psi_1
    =-\sqrt{\mu}\vert D\vert\coth(\sqrt{\mu_2}\D)
    \nabla\psi_1+O(\mu).
    \end{equation}
    In the  BO regime one has $\mu\ll 1$ and $d=0$ (and thus
    $\mu_2=\infty$, $\eps_2=0$), and one gets the approximation
    
    \begin{equation}\label{asympBO}
    {\bf H}^{\mu,d}[\eps\zeta]\psi_1
    \sim-\sqrt{\mu}\vert D\vert
    \nabla\psi_1.
    \end{equation}

The BO regime is thus the limit of the ILW one when the depth of the lower layer is infinite, $d \to 0$, or $\mu_2\to \infty.$

 Using those scalings one derives (see \cite{BLS})  from the full two-layers system the ILW system, written below in horizontal spatial dimension  $N=1,2.$ 

\begin{equation}\label{eqILW}
    \left\lbrace
    \begin{array}{l}
    \dsp[1+\sqrt{\mu}\frac{\alpha}{\gamma}\D\coth(\sqrt{\mu_2}\D)]\dt \zeta+\frac{1}{\gamma}\nabla\cdot ((1-\epsilon\zeta){\bf v})\vspace{1mm}\\
    \indent\dsp-(1-\alpha)\frac{\sqrt{\mu}}{\gamma^2}
    \D\coth(\sqrt{\mu_2}\vert D\vert)\nabla\cdot{\bf v}=0\vspace{1mm},\\
    \dsp \dt {\bf v}
    +(1-\gamma)\nabla\zeta
    -\frac{\epsilon}{2\gamma}\nabla
    \vert{\bf v}\vert^2=0.
    \end{array}\right.
\end{equation}

\begin{remark}
The coefficient $\alpha\geq 0$ in \eqref{eqILW} is a free modeling parameter stemming from the use of the Benjamin-Bona-Mahony (BBM) trick. In dimension $N=1$ and with $\alpha=0$, (\ref{eqILW}) corresponds to
(5.47) of \cite{CGK} which is obtained by expanding the Hamiltonian of the full system with respect to $\epsilon$ and $\mu.$  However this system is not linearly
well-posed.  It is straightforward to ascertain that the condition
$\alpha\geq 1$ insures that (\ref{eqILW}) is linearly well-posed for
either $N = 1$ or $N=2$.
\end{remark}
\begin{remark}
    The ILW equation derived formally in \cite{J,KKD} is obtained as the
    unidirectional limit of the one dimensional
    ($d=1$) version of (\ref{eqILW}) when $\alpha=0, $ see for instance
    \S 5.5 in \cite{CGK}.
\end{remark}

Using the above approximation of the interface operator, leads to the BO system:

\begin{equation}\label{eqBOsyst}
    \left\lbrace
    \begin{array}{l}
    [1+\sqrt{\mu}\frac{\alpha}{\gamma}\D]\dt \zeta+\frac{1}{\gamma}\nabla\cdot ((1-\eps\zeta){\bf v})
    -(1-\alpha)\frac{\sqrt{\mu}}{\gamma^2}
    \D\nabla\cdot{\bf v}=0,\\
    \dt {\bf v}
    +(1-\gamma)\nabla\zeta
    -\frac{\eps}{2\gamma}\nabla
    \vert{\bf v}\vert^2=0
    \end{array}\right.\end{equation}

where $\alpha$ has the same significance as in the previous remark and which again reduces to the BO equation for one-dimensional, unidirectional waves.

The well-posedness of the Cauchy problem for \eqref{eqILW}, \eqref{eqBOsyst} (in space dimension one and two) on long time scales $O(1/\epsilon)$ has been established in \cite{LiXu} together with the limit of solutions of the ILW system to those of the BO system when $\mu_2\to \infty.$  We also refer to \cite{AS} for a study of solitary wave solutions to 
both the one-dimensional ILW and BO systems and to \cite{BDM} for numerical simulations..

\begin{remark}
The above derivation was performed for purely gravity waves. Surface tension effects result in adding a third order dispersive term in the asymptotic models. One gets for instance the so-called Benjamin equation (see \cite{TBB2, TBB3}):

\begin{equation}\label{Benjamin}
u_t+uu_x-Hu_{xx}-\delta u_{xxx}=0,
\end{equation}

where $\delta>0$ measures the capillary effects.

This equation, which is in some sense {\it close} to the KdV equation, is not known to be integrable. Its solitary waves, the existence of which was proven in \cite{TBB2, TBB3} by the degree-theoretic approach, present oscillatory tails. We refer to  \cite{Lin2} for the Cauchy problem in $L^2$ and to \cite{ABR, Ang} for further results on the existence and stability of solitary wave solutions and to \cite{CaAk} for numerical simulations.

In presence of surface tension, the ILW equation has to be modified in the same way. We are not aware of mathematical results on the resulting equation.
\end{remark}

\begin{remark}
The BO equation was fully justified in \cite{OI} as a model of long internal waves in a two-fluid system by taking into account the influence of the surface tension at the interface. The existence time of the full two-fluid system is proportional to the surface tension coefficient (which is very small in real oceanographic systems), making the approximation valid only for very short time.
\end{remark}

\begin{remark}l
Natural generalizations of the ILW and BO equations arise when looking for weak transverse effects, aiming for instance to understand the transverse instability of the BO or ILW solitons and the oblique interactions of such solitary waves, see Section 6.8. In  a weakly transverse regime this leads to Kadomtsev-Petviashvili (KP) versions of the BO and ILW equations (see \cite {AbSe, GZ, Mat12} for the derivation, \cite{LPS2, LPS4, H-GM} for a mathematical study). We will go back in more details to this issue in Section 6.8.

\end{remark}

\begin{remark}
In \cite{Mat17} Matsuno considers the two-fluid system when the upper layer is large with respect to the lower one and in presence of a non trivial topography of the fixed bottom. In this context he derives formally a  forced ILW and a forced BO (fBO) equation and describes the effect on the soliton dynamics in the case of the fBO equation.
\end{remark}
\vspace{0.5cm}
The paper will be organized as follow. The next section will recall the formal framework of Inverse Scattering for both the BO and ILW equations. Section 3 will be devoted on results obtained by pure PDE techniques while Section 4 will describe the rigorous results on the Cauchy problem obtained by IST methods for the BO equation.

Finally in  the last section we  will briefly comment on the periodic case and on related (non integrable) equations and systems, making conjectures, based on numerical simulations on the long time behavior of solutions to ILW, BO and related equations. We will also discuss various issues related to BO and ILW equations : zero dispersion limit, controllability, transverse stability of solitary waves, modified and higher order equations, interaction of solitary waves.

We provide an extensive bibliography since we aim to cover the relevant papers dealing, under various aspects, with the BO and ILW equations.


\vspace{0.5cm}
\noindent{\bf Notations.} We will denote $|\cdot|_p$ the norm in the Lebesgue space $L^p(\R^2),\; 1\leq p\leq \infty$ and $||\cdot||_s$ the norm in the Sobolev space $H^s(\R),\; s\in \R.$ We will use the weighted Sobolev space $H^{s,k}(\R)=\lbrace f\in \mathcal S'(\R), |\langle x\rangle^k\langle i\partial_x\rangle^sf|_2<\infty \rbrace,$  where $\langle x\rangle =(1+x^2)^{1/2}.$ We will denote $\hat {f}$ or $\mathcal F(f)$ the Fourier transform of a tempered distribution $f.$ For any $s\in \R,$ we define $D^s f$  by  its Fourier transform 
$$\widehat{D^s f}(\xi)=|\xi|^s \hat{f}(\xi).$$

\section{An overview of the  Inverse Scattering framework for the ILW and BO equations}

It turns out that the ILW and BO equations are among the relatively few physically relevant equations possessing an inverse scattering formalism.  We give in this Section a brief historical overview and refer to Section 4 for technical details and rigorous and more recent results.

The IST formalism for the BO equation has been given in \cite{AFA, AF}, see also \cite{KaMat, Xu} for extensions to a larger class of potentials. However, before those pioneering works, some facts have been discovered showing the (unexpected) rich structure of the BO equation.

 R.I. Joseph, K.M. Case, A. Nakamura  and  Y. Matsuno seem to have been the first authors in the late seventies to notice the specific properties of the BO equation that lead  to conjecture its complete integrability.  Joseph \cite{J2} derived the expression of the 2-soliton. Case \cite{Ca2} derived also the 2-soliton and in \cite{Ca, Ca3} suggested the existence of an infinite number of conserved quantities $I_n$ and computed the first non trivial ones, namely
 
 \begin{equation}
  I_4=\int\lbrace u^4/4+\frac{3}{2}u^2H(u_x)+2(u_x)^2\rbrace dx,
  \end{equation}
  \begin{equation}
  I_5=\int\lbrace u^5/5+[\frac{4}{3}u^3H(u_x)+u^2H(uu_x)]+[2u(H(u_x))^2+6u(u_x)^2]-4u_{xx}H((u_x)\rbrace dx ,
  \end{equation}
  \begin{equation}
  \begin{split}
  I_6=&\int\lbrace u^6/6+[\frac{5}{4}u^4H(u_x)+\frac{5}{3}u^3H(uu_x)]\\
  &+\frac{5}{2}[5u^2(u_x)^2+u^2(H(u_x)^2+2uH(u_x)H(uu_x)]\\
  &-10[(u_x)^2H(u_x)+2uu_{xx}H(u_x)]+8(u_{xx})^2\rbrace dx.
\end{split}
\end{equation}
This lead Case to conjecture the complete integrability of BO, in particular the existence of a Lax pair.  

We refer to \cite{MX2} for a nice review of the various methods used to derive the BO conservation laws.
\begin{remark}
The "standard" invariants $I_2$ and $I_3$ are respectively given by the $L^2$ norm and the hamiltonian, namely

$$I_2=\frac{1}{2}\int u^2 dx,$$

$$I_3=\int \left(\frac{1}{2}|D^{1/2}u|^2-\frac{1}{6}u^3\right ) dx.$$
\end{remark}


Nakamura \cite{Nak} proved the existence of a B\"{a}cklund transform (see also \cite{BK}) and the existence of an infinite number of conserved quantities. Furthermore  he gave an inverse scattering transform of the BO equation.

In \cite{FoFu}  the Benjamin-Ono  equation is shown to posses two non-local linear operators, which generate its infinitely many commuting symmetries and constants of the motion in involution. These symmetries define the hierarchy of the BO equation, each member of which is a hamiltonian system. The above operators are the nonlocal analogues of the Lenard operator and its adjoint for the Korteweg-de Vries equation, see the discussion in  \cite{SW}.

Further progress were made by Y. Matsuno. In \cite{Mat} he followed Hirota's method (\cite{Hi}) to transform the BO equation into a bilinear form and deduced the explicit expression of the N-soliton. Investigating the asymptotic behavior of the N-soliton as $t\to \infty$ he noticed that, contrary to the case of the KdV equation, no phase-shift appears as the result of collisions of solitons. Still using the bilinear transformation method he gave in \cite{Mat8} the exact formula for the N-soliton solution to the higher order BO equation

$$u_t=\partial_x(\text{grad}\; I_5(u)).$$

In \cite{Mat11} Matsuno  proved, using a recurrence formula derived from the B\"{a}cklund transformation of the BO equation that the functional derivative of a conserved quantity of the BO equation is a conserved quantity. The interaction of N-solitons is studied in more details in \cite{Mat2}. More specifically, the nature of the interaction of 2-solitons is characterized by the amplitudes of the two solitons. This paper contains also a precise analysis of the asymptotics as $t\to \infty$ of solutions to the linearized equation (analogous of the Airy equation for the KdV equation). In \cite{MaKa2} Kaup and Matsuno consider the linearization of the BO equation about the N-soliton solution, proving the linear stability of the N-soliton against infinitesimal perturbations. They also obtain a formal large time asymptotics of solutions to dissipative perturbations of the BO equation. Finally in \cite{Mat9, Mat10} one finds the asymptotic behavior of the number density function of solitons in the small dispersion limit. 

In the survey article \cite{Mat5} a detailed description is given to the interaction process of two algebraic solitons using the pole expansion of the solution, in particular to the effects of small perturbations on the overtaking collision of two BO solitons by employing a direct multisoliton perturbation theory. It is shown that the dynamics of interacting algebraic solitons reveal new aspects which have never been observed in the interaction process of usual solitons expressed in terms of exponential functions.

In addition to \cite{Mat5} we refer to the book \cite{Mat3} for a good description of many "algebraic" aspects of the BO and ILW equations : use of the Hirota bilinear transform method (\cite{Hi}), B\"{a}cklund transforms, multi-solitons, BO and ILW hierarchies...

The Lax pair of the BO equation was derived in \cite{Nak2}, and \cite{BK} while the action-angle variables for the BO equation and their Poisson brackets are determined in \cite{KLM2}. The  direct and inverse transforms for BO was formulated by Fokas Anderson and Ablowitz \cite{AFA,AF}. The direct problem is a differential Riemann-Hilbert problem while the inverse problem consists in solving linear Fredholm equations.
Pure soliton solutions are obtained by solving a linear algebraic system whose coefficients depend linearly on $x, t.$
\vspace{0.3cm}

More precisely, the Lax pair for the BO equation writes, following \cite{AF,AFA}

\begin{equation}
    \label{LaxBO}
    \begin{array}{l}
    i\Phi_x^++\lambda(\Phi^+-\Phi^-)=-u\Phi^+, \\
    i\Phi_t^\pm-2i\lambda \Phi_x^\pm+\Phi_{xx}^\pm-i[u]_x^\pm\Phi^\pm=-\nu \Phi^\pm,
\end{array}
    \end{equation}

where

$$[u]^\pm=\pm\frac{u}{2}+\frac{1}{2i}Hu$$

and where $\lambda$ is a constant interpreted as a spectral parameter, $\nu$ is an arbitrary constant and $[u]^+$ and $[u]^-$ are the boundary values of functions analytic in the upper and lower half complex z-planes respectively.

The core of the analysis of the IST for BO is the linear spectral problem associated to  the first equation in \eqref{LaxBO} which can be interpreted as a differential Riemann-Hilbert problem. This equation  yields unique solutions for $\Phi^+$ and $\Phi^-$ provided one imposes some boundary conditions as $z\to \infty$ say in the upper half plane. The choice made in \cite{AF} is that either $\Phi^+(z,t,\lambda)\to 0$ or $1$ as $z\to\infty,\;\Im z>0.$

One can only consider the $(+)$ functions and thus drop the superscript $^+.$ Let (Jost functions) $M,\overline{M}$ denote the "left" eigenfunctions while $N,\overline{N}$ denote the "right" ones. They are specified by the boundary conditions:

$$M\to 1, \overline{M}\to e^{i\lambda x}\;\text{as}\;x\to -\infty;\quad \overline{N}\to 1\;, N\to e^{i\lambda x}\;\text{as}\; x\to \infty.$$

One then can establish the "scattering equation":

\begin{equation}\label{dirscatBO}
M=\overline{N}+\beta(\lambda,t)\theta(\lambda)N,\quad \beta(\lambda,t)=i\int_{-\infty}^\infty u(y,t)M(y,t,\lambda)e^{-i\lambda y} dy,
\end{equation}

where $\theta(y)=1, \lambda >0$ and $0, \lambda <0.$

The evolution of the scattering data is given by 

$$\rho(k;t)=ik^2\exp(ik^2t),$$

$$f(k,t)=f(k,0)\exp(-ik^2t).$$

The solution of the "inverse" problem consists essentially in solving \eqref{dirscatBO}.  The evolution of the scattering data is simple and the solution of the inverse problem is characterized by a linear Fredholm equation.

The soliton solutions are obtained by taking

$$\rho(k;t)=0,\quad N(x,t;k)=1-i\sum_{j=1}^n\frac{\phi_j}{k-k_j}.$$



\vspace{0.5cm}
 The IST formalism for ILW  has been given in \cite{KAS, KSA} but, as in the case of the BO equation, facts suggesting the (formal) complete  integrability of the equation have been discovered before, for instance the existence of multi-solitons (see \cite{JE, NM}).

The ILW equation possesses an infinite sequence of conserved quantities (see {\it eg} \cite{LR, Mat3}) which leads to a ILW hierarchy. The first non trivial one is
\begin{equation}\begin{split} 
J_4(u)=   &\int_{-\infty}^\infty (\frac{1}{4}u^4+\frac{3}{2}u^2\mathcal T(u_x)+\frac{1}{2}u_x^2+\frac{3}{2}[\mathcal T(u_x)]^2\\
&+\frac{1}{\delta}[\frac{3}{2}u^3+\frac{9}{2}u\mathcal T(u_x)]+\frac{3}{2\delta^2})dx.
\end{split}
\end{equation}

The direct 
 scattering problem is associated with a Riemann-Hilbert problem in a strip of the complex plane.

More precisely, the Lax pair for the ILW equation is given (see \cite{KAS, KSA} and Chapter 4 in \cite{AC}) :

 \begin{equation}
    \label{LaxILW}
    \begin{array}{l}
    iv_x^++(u-\lambda)v^+=\mu v^-, \\
    iv_t^++i(2\lambda+\delta^{-1})v_x^\pm +v_{xx}^\pm +(\pm iu_x-\mathcal T(u_x)+\nu)v^\pm=0,
\end{array}
    \end{equation}
    
    where $\lambda, \mu$ are parametrize as
    

    $$\lambda(k)=-\frac{1}{2}k\coth(k\delta)\quad \text{and} \quad \mu(k)=\frac{1}{2}k\;\text{cosech}(k\delta).$$
    
    $k$ is a constant which is interpreted as a spectral parameter and $\nu$ is an arbitrary constant. Given $u$, the first equation in \eqref{LaxILW} defines a Riemann-Hilbert problem in the horizontal  strip, more precisely, $v^\pm(x)$ represent the boundary values of analytic functions in  the strip between $\Im(z)=0$ and $\Im(z)=2\delta, z=x+iy $ and periodically extended in the vertical direction.
    
    Using the operator $\mathcal T, v^\pm$ may be written as 
    
    \begin{equation}
    \label{newILW}
    \begin{array}{l}
    v^+(x)=\lim_{y\to 0} v(x)=\frac{1}{2}(I-i\mathcal T)\psi(x), \\
    v^-(x)=\lim_{y\to 0}v(x)=\frac{1}{2}(I+i\mathcal T)\psi(x),
\end{array}
    \end{equation}
    
    where $\psi\in L^1(\R),$ is H\"{o}lder of exponent $\alpha$ and for $\Im z \neq 0$ (mod $2\delta$) , $v(z)$ is given by

    $$v(z)=\frac{1}{4i\delta}\int_{-\infty}^\infty \coth\lbrace \frac{\pi}{2\delta}(y-z)\rbrace \psi(x)dy.$$
    
    Defining

    $$W^{\pm}(x;k)=v^\pm(x;k)\exp \lbrace\frac{1}{2}ik(x\mp i\delta)\rbrace,$$
    
    \eqref{newILW} can be rewritten as

    \begin{equation}
    \label{newnewILW}
    \begin{array}{l}
    iW_x^++[\zeta_++1/(2\delta)](W^+-W^-)=-uW^+, \\
    iW_t^\pm-2i\zeta_+W_x^\pm+W_{xx}^\pm+[\pm iu_x-\mathcal T(u_x))+\rho]W^\pm=0,
\end{array}
    \end{equation}
    
    where 
    
    $$\zeta_\pm(k)=\frac{1}{2}k\pm\frac{1}{2}k\coth(k\delta)\mp1/(2\delta),$$
    
    and
    
    $$\rho=-k\zeta_++\frac{1}{4}k^2+\nu.$$
    
    The solution of \eqref {newnewILW} is given by the integral equation
    
    \begin{equation}
    W^+(x;k)=W_0(x;k)+\int_{-\infty}^\infty G(x,y;k)u(y)W^+(y;k)dy
    \end{equation}

where $W_0(x;k)$ is the solution of \eqref{newnewILW} corresponding  to $u=0$ and

$G$ is the Green function satisfying

\begin{equation}
i\frac{\partial}{\partial x}\lbrace G^+(x,y;k)\rbrace +[\zeta_++1/(2\delta)][G^+(x,y;k)-G^-(x,y;k)]=-\delta(x-y),
\end{equation}

where $G^\pm(x,y;k)=G(x\mp i\delta,y;k).$

We refer to \cite{AC} for an integral representation of $G_\pm.$ As in the case of the BO equation one needs also the eigenfunction $W^+(x;k).$ We also denote $M,\overline{M}$ the "left" eigenfunctions and $N,\overline{N}$ the "right" eigenfunctions, that have the following asymptotic behavior :

$$M(x;k)\sim 1,\quad \overline{M}(x:k)\sim e^{ikx+k\delta},\quad \text{as}\; x\to -\infty$$

$$N(x;k)\sim e^{ikx+k\delta},\quad \overline{N}(x,k)\sim 1,\quad \text{as}\;x\to \infty.$$

The eigenfunctions $M,N,\overline{N}$ are related through the completeness relation

\begin{equation}\label{compl}
M(x;k)=a(k)\overline{N}(x;k)+b(k)N(x;k),
\end{equation}

where $a(k)$ and $b(k)$ have the integral representation:
\begin{equation}
\begin{array}{l}
a(k)=1+\frac{1}{2i\zeta_+(k)}\int_{-\infty}^\infty u(y)M(y;k) dy,\\
b(k)=-\frac{1}{2i\zeta_+(k)}\theta(\zeta_++1/(2\delta))\int_{-\infty}^\infty u(y)M(y;k)e^{-iky-k\delta}dy.
\end{array}
\end{equation}

\vspace{0.3cm}

One can prove that the evolution of the scattering data is simple, namely

$$a(k,t)= a(k,0) ,\quad b(k,t)=b(k,0)\exp\lbrace ik[k\coth(k\delta)-1/\delta]t\rbrace,$$

so that

$$\rho(k,t)=\rho(k,0)\exp\lbrace ik[k\coth(k\delta)-1/\delta]\rbrace.$$

Similarly,

$$k_j(t)=k_j(0),$$

$$C_j(t)=C_j(0)\exp\lbrace ik_j[k_j\coth(k_j \delta)-1/\delta]t\rbrace,$$

for $j=1,2,...,n.$

The inverse scattering problem is based on equation \eqref{compl}  and consist in reconstructing $M,\ N, \overline{N}$ from the knowledge of $a(k), b(k)$ together with appropriate information about the ground states.  The formal solution of this Riemann-Hilbert problem is given in \cite{AC}.

\begin{remark}
The pure soliton solutions are recovered by taking $\rho(k;t)=0.$
\end{remark}

All the results in this subsection are formal. In particular Fokas and Ablowitz \cite{AF} have to make some key spectral assumptions in their definition of the scattering data of the IST for the BO equation. Also, their IST framework does not behave well enough to be solved by iteration.

We refer to Section 4 for a review on rigorous results for the Cauchy problem and related issues using IST techniques.

\section{Rigorous results by PDE methods}

\subsection{The linear group}
We will recall here the various dispersive estimates satisfied by the linear BO and ILW equations.  The linear part of both the BO and the ILW equation defines a unitary group $S_{BO}(t),$ resp. $S_{ILW}(t)$ in all Sobolev spaces $H^s(\R), s\geq 0$ which is by Plancherel unitarily equivalent in $L^2$ to the multiplication by, respectively,  $e^{it\xi|\xi|}$  and $e^{it(\xi^2\coth(\delta \xi)-\xi/\delta)}.$

The fundamental solutions:

$$G_{BO}(x,t)=\mathcal F^{-1}(e^{it\xi|\xi|})(x)\quad \text{ and}\quad  G_{ILW}(x,t)=\mathcal F^{-1}(e^{it(\xi^2\coth(\delta \xi)-\xi/\delta)})(x),$$

play an important role in the dispersive properties of the BO and ILW equations.  $G_{BO}(x,1)$ is easily seen to be a bounded $C^\infty$ function. Its asymptotic behavior is obtained in \cite{SSS} (particular case of Theorem 3.1):

$$G_{BO}(x,1)\sim \sqrt \pi\cos\left(\left(\frac{|x|}{2}\right)^2+\frac{\pi}{4}\right),\;\text{as}\; x\to -\infty,$$

and  

$$G_{BO}(x,1)\sim \sum_{k=0}^\infty(-1)^k\frac{[2(2k+1)]!}{ (2k+1)!x^{2(2k+1)+1}}\;\text{as}\; x\to +\infty.$$



On the other hand one proves that:
\begin{equation} 
 |G_{ILW}(x,t)| \lesssim
\begin{cases}
  t^{-\frac{1}{3}} \left< t^{-1/3}x\right>^{-1/4} \text{ for } x\leq t \\
 t^{-\frac{1}{2}} \text{ for } x\geq t .
  \end{cases}
\end{equation}


In both cases one obtains  the dispersive estimate 

\begin{equation}\label{disp}
|S_{BO}(t)\phi|_\infty, \;|S_{ILW}(t)\phi|_\infty\lesssim \frac{1}{t^{1/2}}|\phi|_1
\end{equation}

yielding  for instance the same Strichartz estimates as the linear one-dimensional Schr\"{o}dinger group.

The BO and ILW groups display a Kato type local smoothing   property \cite{CoSa, KPV}. The optimal results for the group $S_{BO}(t)$ are gathered in the next theorem from \cite{KPV2}:

\begin{theorem}\label{KPV}
There exist constants $c_0, c_1$ such that
\begin{equation}
\left(\int_{-\infty}^\infty|D^{1/2}S_{BO}(t)u_0(x)dt\right)^{1/2}=c_0|u_0|_2.
\end{equation}
for any $x\in \R.$
\begin{equation}
|D^{1/2}\int_{-\infty}^\infty S_{BO}(t)g(.,t)dt|_2\leq c_1\int_{-\infty}^\infty\left(\int_{-\infty}^\infty |g(x,t)|^2dt\right)^{1/2}dx,
\end{equation}
and
\begin{multline}
\sup_x\left(\int_{-\infty}^\infty|\partial_x\left(\int_0^tS_{BO}(t-\tau)f(.,\tau)d\tau\right)|^2dt\right)^{1/2}\\
\leq c_1\int_{-\infty}^\infty\left(\int_{-\infty}^\infty|f(x,t)|^2dt\right)^{1/2}dx.
\end{multline}

\end{theorem}

As aforementioned the Strichartz estimates are consequences of the dispersive estimate \eqref{disp}. We state here the version in \cite{KPV2}
\begin{theorem}\label{Stri}
Let $p\in [2,\infty)$ and q be such that $2/q=1/2-1/p.$ Then
\begin{equation}\label{dips1}
\left(\int_{-\infty}^\infty|S_{BO}(t)u_0|_p^q dt\right)^{1/q}\leq c|u_0|_2
\end{equation}
and 
\begin{equation}\label{Sri2}
\left(\int_{-\infty}^\infty |\int_0^tS_{BO}(t-\tau)f(.,\tau)d\tau |_p^q\right)^{1/q}\leq C\left(\int_{-\infty}^\infty|f(.,\tau)|_{p'}^{q'}dt\right)^{1/q'}
\end{equation}

where $1/p+1/p'=1/q+1/q'=1$.
\end{theorem}

The last group of dispersive estimates  are those on the maximal function $\sup_t S_{BO}(t).$ Then, see \cite{KPV2}:

\begin{theorem}\label{max}
\begin{equation}\label{max1}
\left( \int_{-\infty}^\infty \sup_{-\infty<t<\infty}|S_{BO}(t)u_0(x)|^4dx\right)^{1/4}\leq c|D^{1/4}u_0|_2,
\end{equation}

and
\begin{multline}\label{max2}
\left( \int_{-\infty}^\infty \sup_{-\infty<t<\infty}|\int_0^t S_{BO}(t-\tau)f(.,\tau)d\tau|^4dx\right)^{1/4}\\
\leq c\left(\int_{-\infty}^\infty\left(\int_{-\infty}^\infty |D^{1/2}f(x,t)|dt\right)^{4/3}dx\right)^{3/4}.
\end{multline}
Moreover for $s>\frac{1}{2}$ and $\rho>\frac{3}{4}$
\begin{equation}\label{max3}
\left(\int_{-\infty}^\infty \sup_{0\leq t\leq T}|S_{BO}(t)u_0(x)|^2\right)^{1/2}\leq c(1+T)^\rho||u_0||_s.
\end{equation}

\end{theorem}
\begin{remark}
As will be discussed below, and contrary to the case of the KdV equation, the above estimates cannot be used to define a functional space on which one could implement an iterative scheme based on the Duhamel formulation
$$u(t)=S_{BO}(t)u_0+\int_0^t S_{BO}(t-\tau)uu_x(\tau)d\tau.$$
 This is due to the quasilinear nature of the BO equation.
\end{remark}

\subsection{An easy result}
Being skew-adjoint perturbations of the inviscid Burgers equation, the ILW and BO equations are easily seen to be locally well-posed in $H^s(\R), s>\frac{3}{2},$  (\cite{Io, JCS}). It turns out that  they are actually {\it globally} well posed in the same range of Sobolev spaces.

The argument goes back to \cite{JCS} and was applied in this context in \cite{ABFS}. We describe it briefly for the BO equation. The following computations are formal but can be easily justified by smoothing the equation and/or the initial data. We start from the identity 

\begin{equation}\label{est1}
\frac{1}{2}\frac{d}{dt}|D^su|^2_0+(D^s(uu_x),D^su)=0.
\end{equation}

For $s>\frac{3}{2}$ one deduces after some manipulations using commutator estimates and the Sobolev imbedding

$$|u_x|_\infty\leq \frac{c}{\sqrt \eta}||u||_{3/2+\eta},\; \forall \eta>0,$$

that

$$|(D^s(uu_x),D^su)|\leq \frac{c}{\sqrt \eta}||u||_{3/2+\eta}||u||^2_s,$$

holding for any $\eta>0$ such that $\frac{3}{2}+\eta<s.$  For such an $\eta$ we have the interpolation inequality 

$$||u||_{3/2+\eta}\leq c||u||_s^{2\gamma \eta}||u||_{3/2}^{1-2\gamma \eta}$$

where $1-\theta=\frac{2\eta}{2s-3}=:2\gamma \eta$ and c is a constant independent of $\eta.$

On the other hand, using the conservation law $I_5$ one can establish the uniform bound

\begin{equation}\label{est3}
||u||_{L^\infty (\R;H^{3/2}(\R))}\leq ca(||u_0||_{3/2}).
\end{equation}

Combining those estimates yields

\begin{equation}\label{est4}
\frac{1}{2}\frac{d}{dt}||u||_s^2\leq C\left(\frac{1}{\sqrt \eta}||u||_s^{2+2\gamma \eta}\right),
\end{equation}

where the constant $C$ is defined by $C=[ca(||u_0||_{3/2})]^{1-2\gamma \eta}.$

Integrating this last inequality in time  leads to the estimate

\begin{equation}\label{est5}
||u(\cdot,t)||_s^2\leq y(t),
\end{equation}

where $y$ is the solution of the differential equation

$$y'(t)=\frac{C}{\sqrt \eta}y(t)^{1+\gamma \eta},\; y(0)=||u_0||^2_s$$

on its maximal interval of existence $[0, T(\eta)].$  Here $\gamma=1/(2s-3).$ This equation is easily integrated, finding that

$$y(t)=(||u_0||_s^{-2\gamma \eta}-\gamma\sqrt \eta Ct)^{-1/\gamma \eta},$$

whence

$$T(\eta)=\frac{1}{\gamma C\sqrt \eta}||u_0||_s^{-2\gamma\eta}\to \infty \quad\text{as}\;\eta\to 0.$$

For any fixed $T>0$ we can choose $\eta>0$ so small that $T<\frac{1}{2}T(\eta).$ Then it follows that for  $0\leq t\leq T,$

$$y(t)\leq c(T;||u_0||_{3/2})||u_0||_s^2.$$

This implies an a priori bound that is crucial to prove the existence of a unique global solution $u\in L^\infty_{\text{loc}}(\R, H^s(\R))$ emerging from an initial data $u_0\in H^s(\R), s>3/2.$

The strong continuity in time and the continuity of the flow map is established by using the Bona-Smith method \cite{BoSm}.
\begin{remark}
1. A similar global result holds true for the ILW equation, see \cite{ABFS}.

2. As aforementioned, it is  proved in \cite{ABFS} that the $H^s, s>3/2$ solutions of the ILW equation converge to that of the BO (resp. KdV) equation when $\delta\to +\infty$ (resp. $\delta \to 0$). Recall that in the KdV case, (and this is often missed in the literature...) one has to rescale the ILW solution $u_\delta$ as 

$$v_\delta(x,t)=\frac{3}{\delta}u_\delta(x,\frac{3}{\delta}t),$$

and then let $\delta \to 0,$ to obtain a KdV solution.
\end{remark}
\subsection{Global weak solutions}
By deriving  local smoothing properties {\it \`a la Kato} for the nonlinear equation and various delicate commutator estimates, Ginibre and Velo (\cite{GV,GV2,GV3, GV4})  proved various global existence results of weak solutions for a class of generalized BO equations. For the BO equation itself (and also for the ILW  equation \footnote{the results hold also for a rather general class of nonlocal dispersive equations, see \cite{GV4}.}), this implies the existence of a global weak solution in $L^\infty(\R; L^2(\R))\cap L^2_{\text{loc}}((\R; H^{1/2}_{\text{loc}}(\R))$ for any initial data in $L^2(\R)$ and similar results for data with higher regularity, {\it eg} $H^1(\R)$.Those results hold also true for the ILW equation (see \cite{GV4} Section 6).

For initial data in the energy space $H^{1/2}(\R)$  the existence of  global weak solutions in the space  $L^\infty(\R; H^{1/2}(\R))$ has been established in \cite{JCS} for both the BO and ILW equations.

\subsection{Semilinear versus quasilinear}

Deciding whether a nonlinear dispersive equation is semilinearly well-posed is somewhat subtle. The distinction plays a fundamental role in the choice of the method for solving the Cauchy problem.

To start with we consider the case of the Benjamin-Ono equation following \cite{MST4}. Since they are relatively simple we will give complete proofs.

We thus consider the Cauchy problem
\begin{equation}\label{BObis}
\left \{
\begin{array}{l}
u_{t}-Hu_{xx}+uu_{x} =  0, ~(t,x)\in\R^2,\\
u(0,x)  =  \phi(x).
\end{array}
\right.
\end{equation}
Setting $S(t)= e^{itH\partial^2_x},$ we write (\ref{BObis}) as an integral equation:
\begin{eqnarray}\label{ie}
u(t)=S(t)\phi-\int_{0}^{t}S(t-t')(u_{x}(t')u(t'))dt'.
\end{eqnarray}

The main result  is the following

\begin{theorem}\label{th1}
Let $s\in\R$  and $T$ be a positive real number.
Then there does not exist a space $X_T$ continuously
embedded in $C([-T,T],H^{s}(\R))$
such that there exists $C>0$ with
\begin{eqnarray}\label{contr1}
\|S(t)\phi\|_{X_T}\leq C\|\phi\|_{H^{s}(\R)},\quad \phi\in
H^{s}(\R),
\end{eqnarray}
and
\begin{eqnarray}\label{contr2}
\left\|\int_{0}^{t}S(t-t')\left[u(t')u_{x}(t')\right]dt'\right\|_{X_T}
\leq
C\|u\|_{X_T}^{2},\quad u\in X_T.
\end{eqnarray}
\end{theorem}

Note that (\ref{contr1}) and (\ref{contr2}) would be needed to
implement a Picard iterative scheme (actually the second iteration) on (\ref{ie}), in the space $
X_T $. As a consequence of Theorem \ref{th1} we can obtain the
following result.

\begin{theorem}\label{th2}
Fix $s\in\R$. Then there does not exist a $T>0$ such that
{\rm{(\ref{BObis})}} admits a unique local solution defined on the
interval $[-T,T]$ and such that the flow-map data-solution $$
\phi\longmapsto u(t), \quad t\in [-T,T],$$ for
{\rm{(\ref{BObis})}} is $C^{2}$ differentiable at zero from
$H^{s}(\R)$ to $H^{s}(\R)$.
\end{theorem}

\begin{remark}
This result implies that the Benjamin-Ono equation is "quasilinear". It has been precised in \cite{KT2} where it is shown that the flow map cannot even be locally Lipschitz in $H^s$ for $s\geq 0.$ This is of course in strong contrast with the KdV equation.
\end{remark}

\subsubsection{Proof of Theorem \ref{th1}}\label{sec2.1}
Suppose that there exists a space $X_T$ such that (\ref{contr1})
and (\ref{contr2}) hold. Take $u=S(t)\phi$ in (\ref{contr2}). Then
$$ \left\|\int_{0}^{t}S(t-t')\left[(S(t')\phi)
(S(t')\phi_{x})\right]dt'\right\|_{X_T} \leq
C\|S(t)\phi\|_{X_T}^{2}. $$ Now using (\ref{contr1}) and that
$X_T$ is continuously embedded in $C([-T,T],H^{s}(\R))$ we obtain
for any $t\in[-T,T]$ that
\begin{eqnarray}\label{fail}
\left\|\int_{0}^{t}S(t-t')\left[(S(t')\phi)
(S(t')\phi_{x})\right]dt'\right\|_{H^{s}(\R)}
\lesssim
\|\phi\|_{H^s(\R)}^{2}.
\end{eqnarray}
We show that (\ref{fail}) fails by choosing an appropriate $\phi$.

Take $\phi$ defined by its Fourier
transform as\footnote{The analysis below works as well for
${\mathcal Re}\, \phi$ instead of $\phi$ (some new harmless terms appear).}
$$
\widehat{\phi}(\xi)
=
\alpha^{-\frac{1}{2}} {\bf 1}_{I_1}(\xi) + \alpha^{-\frac{1}{2}}N^{-s}
{\bf 1}_{I_2}(\xi),\quad N\gg 1,~~ 0<\alpha\ll 1, $$ where $I_1$,
$I_2$ are the intervals $$ I_{1}=[\alpha/2,\alpha],~~ I_{2}=
[N,N+\alpha]. $$ Note that $\|\phi\|_{H^{s}}\sim 1$. We will use
the next lemma.

\begin{lemma}\label{ff}
The following identity holds:
\begin{multline}
\int_{0}^{t}S(t-t')\Bigl[ (S(t') \phi) (S(t') \phi_x) \Bigr] \,
dt' =\\ \int_{\R^2} e^{ix\xi+it p(\xi) } \, \xi \,
\hat{\phi}(\xi_1) \hat{\phi}(\xi-\xi_1) \frac{e^{it(p(\xi_1) +
p(\xi-\xi_1) -p(\xi))}-1} {p(\xi_1)+p(\xi-\xi_1) -p(\xi)} d\xi
d\xi_{1},
\end{multline}
where $ p(\xi)=\xi |\xi| $.
\end{lemma}

{\it Proof of Lemma} \ref{ff}. 
Taking the inverse Fourier transform with respect to $ x $, it is
easily seen
 that

According to the above lemma,
$$
\int_{0}^{t}S(t-t')\left[(S(t')\phi)
(S(t')\phi_{x})\right]dt'
=
c(f_{1}(t,x)+f_{2}(t,x)+f_{3}(t,x)), $$ where, from the definition
of $ \phi $,  we have the following representations for $f_1$,
$f_2$, $f_3$:
\begin{eqnarray*}
f_{1}(t,x)
& = &
\frac{c}{\alpha}
\int_{{\xi_{1}\in I_{1}}\atop
{\xi-\xi_{1}\in I_{1}}}
\xi\,
e^{ix\xi+it\xi|\xi|}
\frac{e^{it(\xi_1|\xi_1|+(\xi-\xi_1)|\xi-\xi_1|-\xi|\xi|)}-1}
{\xi_1|\xi_1|+(\xi-\xi_1)|\xi-\xi_1|-\xi|\xi|}
d\xi d\xi_{1},
\\[0.5cm]
f_{2}(t,x)
& = &
\frac{c}{\alpha N^{2s}}
\int_{{\xi_{1}\in I_{2}}\atop
{\xi-\xi_{1}\in I_{2}}}
\xi\,
e^{ix\xi+it\xi|\xi|}
\frac{e^{it(\xi_1|\xi_1|+(\xi-\xi_1)|\xi-\xi_1|-\xi|\xi|)}-1}
{\xi_1|\xi_1|+(\xi-\xi_1)|\xi-\xi_1|-\xi|\xi|}
d\xi d\xi_{1},
\\[0.5cm]
f_{3}(t,x)
& = &
\frac{c}{\alpha N^{s}}
\int_{{\xi_{1}\in I_{1}}\atop
{\xi-\xi_{1}\in I_{2}}}
\xi\,
e^{ix\xi+it\xi|\xi|}
\frac{e^{it(\xi_1|\xi_1|+(\xi-\xi_1)|\xi-\xi_1|-\xi|\xi|)}-1}
{\xi_1|\xi_1|+(\xi-\xi_1)|\xi-\xi_1|-\xi|\xi|}
d\xi d\xi_{1}
\\[0.5cm]
& &
+~~~
\frac{c}{\alpha N^{s}}
\int_{{\xi_{1}\in I_{2}}\atop
{\xi-\xi_{1}\in I_{1}}}
\xi\,
e^{ix\xi+it\xi|\xi|}
\frac{e^{it(\xi_1|\xi_1|+(\xi-\xi_1)|\xi-\xi_1|-\xi|\xi|)}-1}
{\xi_1|\xi_1|+(\xi-\xi_1)|\xi-\xi_1|-\xi|\xi|}
d\xi d\xi_{1}.
\end{eqnarray*}
Set
$$
\chi(\xi,\xi_1):=\xi_1|\xi_1|+(\xi-\xi_1)|\xi-\xi_1|-\xi|\xi|.
$$
Then clearly
\begin{eqnarray*}
{\mathcal F}_{x\mapsto \xi}(f_{1})(t,\xi)
& = &
\frac{c\,\xi e^{it\xi|\xi|}}
{\alpha}
\int_{{\xi_{1}\in I_{1}}\atop
{\xi-\xi_{1}\in I_{1}}}
\frac
{e^{it\chi(\xi,\xi_1)}-1}
{\chi(\xi,\xi_1)}
d\xi_{1},
\\[0.5cm]
{\mathcal F}_{x\mapsto \xi}(f_{2})(t,\xi)
& = &
\frac{c\,\xi e^{it\xi|\xi|}}
{\alpha N^{2s}}
\int_{{\xi_{1}\in I_{2}}\atop
{\xi-\xi_{1}\in I_{2}}}
\frac
{e^{it\chi(\xi,\xi_1)}-1}
{\chi(\xi,\xi_1)}
d\xi_{1},
\\[0.5cm]
{\mathcal F}_{x\mapsto \xi}(f_{3})(t,\xi)
& = &
\frac{c\,\xi e^{it\xi|\xi|}}
{\alpha N^{s}}
\left(
\int_{{\xi_{1}\in I_{1}}\atop
{\xi-\xi_{1}\in I_{2}}}
\frac
{e^{it\chi(\xi,\xi_1)}-1}
{\chi(\xi,\xi_1)}
d\xi_{1}
+
\int_{{\xi_{1}\in I_{2}}\atop
{\xi-\xi_{1}\in I_{1}}}
\frac
{e^{it\chi(\xi,\xi_1)}-1}
{\chi(\xi,\xi_1)}
d\xi_{1}\right).
\end{eqnarray*}
Since the supports of ${\mathcal F}_{x\mapsto \xi}(f_{j})(t,\xi)$,
$j=1,2,3$, are disjoint, we have $$
\left\|\int_{0}^{t}S(t-t')\left[(S(t')\phi)
(S(t')\phi_{x})\right]dt'\right\|_{H^{s}(\R)} \geq
\|f_{3}(t,\cdot)\|_{H^{s}(\R)}. $$ We now give a lower bound for
$\|f_{3}(t,\cdot)\|_{H^{s}(\R)}$. Note that for
$(\xi_1,\xi-\xi_1)\in I_1\times I_2$ or $(\xi_1,\xi-\xi_1)\in
I_2\times I_1$ one has $|\chi(\xi,\xi_1)|=2|\xi_1(\xi-\xi_1)|\sim
\alpha N$. Hence it is natural to choose $\alpha$ and $N$ so that
$\alpha N= N^{-\epsilon}$, $ 0 < \epsilon\ll 1 $. Then $$ \left|
\frac {e^{it\chi(\xi,\xi_1)}-1} {\chi(\xi,\xi_1)} \right|=
|t|+O(N^{-\epsilon}) $$ for $\xi_{1}\in I_{1}$, $\xi-\xi_{1}\in
I_{2}$ or $\xi_{1}\in I_{2}$, $\xi-\xi_{1} \in I_{1}$. Hence for
$t\neq 0$,
\begin{eqnarray*}
\|f_{3}(t,\cdot)\|_{H^{s}(\R)}
\gtrsim
\frac{N\, N^{s}\,\alpha\,\alpha^{\frac{1}{2}}}
{\alpha N^{s}}
=
\alpha^{\frac{1}{2}}N.
\end{eqnarray*}
Therefore we arrive at
\begin{eqnarray*}
1\sim \|\phi\|^{2}_{H^{s}(\R)}
\geq
\|f_{3}(t,\cdot)\|_{H^{s}(\R)}
 \geq
\alpha^{\frac{1}{2}}N\sim N^{\frac{1-\epsilon}{2}},
\end{eqnarray*}
which is a contradiction for $N\gg 1$ and $\epsilon\ll 1$. This
completes the proof of Theorem \ref{th1}. \qquad
\endproof

\subsection{Proof of Theorem \ref{th2}}

Consider the Cauchy problem
\begin{equation}\label{BO}
\left \{
\begin{array}{l}
u_{t}-Hu_{xx}+uu_{x} =  0, \\
u(0,x)  =  \gamma \phi,~~ \gamma\ll 1, \quad \phi \in H^s(\R) \;.
\end{array}
\right.
\end{equation}
Suppose that $u(\gamma,t,x)$ is a local solution of (\ref{BO})
and that the flow map is $C^2$ at the origin from $H^{s}(\R)$
to $H^{s}(\R)$. We have successively

$$u(\gamma, t,x)=\gamma S(t)\phi+\int_0^t S(t-t') u(\gamma,t',x)u_x(\gamma,t',x)dt'$$

$$\frac{ \partial u }{\partial \gamma} (0,t,x)=S(t)\phi (x)=:u_1(t,x)$$
$$
\frac
{\partial^{2}u}{\partial\gamma^{2}}(0,t,x)
=
-2\int_{0}^{t}S(t-t')\left[(S(t')\phi) (S(t')\phi_{x})\right]dt'.
$$ The assumption of $C^2$ regularity yields $$
\left\|\int_{0}^{t}S(t-t')\left[(S(t')\phi)
(S(t')\phi_{x})\right]dt'\right\|_{H^{s}(\R)} \lesssim
\|\phi\|_{H^s(\R)}^{2}. $$ But the above estimate is (\ref{fail}),
which has been shown to fail in section \ref{sec2.1}. \qquad
\endproof

\begin{remark}
The previous results are in fact valid in a more general context.
We consider now the class of equations
\begin{eqnarray}\label{general}
u_t+uu_x-Lu_x =  0,\quad u(0,x) & = & \phi(x),\quad (t,x)\in\R^2,
\end{eqnarray}
where $L$ is defined via the Fourier transform
$$
\widehat{Lf}(\xi)
=
\omega(\xi)\hat{f}(\xi).
$$
Here $\omega(\xi)$ is a continuous real-valued function.
Set $p(\xi)=\xi\,\omega(\xi)$. We assume that $p(\xi)$ is
differentiable and such that, for some $ \gamma \in \R $,
\begin{eqnarray}\label{borne}
|p'(\xi)|
\lesssim
|\xi|^{\gamma},\quad \xi \in \R.
\end{eqnarray}

The next theorem shows that (\ref{general}) shares the bad
behavior of the Benjamin--Ono equation with respect to iterative
methods.

\begin{theorem} \label{g}
Assume that {\rm{(\ref{borne})}} holds with $\gamma\in [0,2[$.
Then the conclusions of Theorems {\rm{\ref{th1}}}, {\rm{\ref{th2}}} are valid for the
Cauchy problem {\rm{(\ref{general})}}.
\end{theorem}

The proof follows the considerations of the previous section.
The main point in the analysis is that for $\xi_1\in I_1$,
$\xi-\xi_1\in I_2$ one has
$$
|p(\xi_1)+p(\xi-\xi_1)-p(\xi)|\lesssim \alpha N^{\gamma}, \quad \quad
 \alpha \ll 1, \quad N\gg 1.
$$ We choose $\alpha$ and $N$ such that $\alpha N^{\gamma}=
N^{-\epsilon}$, $0<\epsilon\ll 1$. We take the same $\phi$ as in
the proof of Theorem \ref{th1} and arrive at the lower bound $$
1\sim \|\phi\|^{2}_{H^{s}(\R)} \geq \alpha^{\frac{1}{2}}N=
N^{1-\frac{\gamma+\epsilon}{2}},$$ which fails for $0<\epsilon\ll
1,\,\gamma\in [0,2[$.

Here we give several examples where Theorem \ref{g} applies.

$\bullet$ Pure power dispersion: $$ \omega(\xi)=|\xi|^{\gamma},
\quad 0\leq \gamma <2. $$ This dispersion corresponds to a class
of models for vorticity waves in the
 coastal zone (see \cite{Sm}).
It is interesting to notice that the case $\gamma=2$ corresponds
to the KdV equation which can be solved by iterative methods, see {\it eg} 
\cite{KPV3}. Therefore Theorem \ref{g} is sharp for a pure
power dispersion. However, the Cauchy problem corresponding to
$1\leq \gamma <2$ has been proven in \cite{HIKK} to
be locally well-posed by a compactness method combined with sharp
estimates on the linear group for initial data in $H^{s}(\R)$,
$s\geq (9-3\gamma)/4$.

$\bullet$ Perturbations of the Benjamin-Ono equation:

$\omega(\xi) =
(|\xi|^2+1)^{\frac{1}{2}}.$ This case corresponds to an equation introduced by Smith \cite{Sm} for continental shelf waves.

$\omega(\xi) =
\xi\, \text{coth}(\xi)$. This corresponds to the Intermediate long wave
equation, cf. \cite{KKD, Liu, ABS, ABFS} and Chapter 5 Section 2.

\end{remark}

\begin{remark}
As it is clear from the above proof, the "ill-posedness" of the Benjamin -Ono equation is due to the "bad" interactions of very small and very large frequencies. This phenomena does not occur in the periodic case, for initial data say of zero mean, see  Subsection 5.2 below.

This is in contrast with similar "ill-posedness" results for the KP-I equation (see \cite{MST3}) which are due to the large zero set of a resonant function and which persist in the periodic case.
\end{remark}
\begin{remark}
The generalized BO equation
\begin{equation}\label{geneBO}
u_t+u^pu_x-Hu_{xx}=0, \; p\geq 2,
\end{equation}

is in fact {\it semilinear} since the Cauchy problem for small data in suitable Sobolev spaces was proven in \cite{KPV2} to be locally-well posed by a contraction method.
\end{remark}

\subsection{Global well-posedness in $L^2$}
The global well-posedness result in $H^s, s>3/2$ does not use any dispersive estimate (but it uses the existence of a non trivial invariant). The results of Ginibre and Velo use a Kato type (dispersive) local smoothing  but they concern only weak solutions. The next step was to use in a more crucial way the dispersion properties to obtain the local well-posedness (LWP) of the Cauchy problem in larger Sobolev spaces, aiming to reach at least the energy space $H^{1/2}(\R).$ 

The first significant result in that direction was from Koch and Tzvetkov \cite{KT} who proved the LWP in $H^s(\R), s>5/4.$ \footnote{Ponce \cite{Po} used dispersive properties but only reached $H^{3/2}(\R).$}  The main idea of the proof is to improve the dispersive estimates (Strichartz estimates) by localizing them in space frequency dependent time intervals together with classical energy estimates. This method was improved by Kenig and Koenig \cite{KeKo} who obtained LWP in the space $H^s(\R), s>9/8.$


A breakthrough was made by T. Tao \cite{Tao2} who obtained the LWP (and thus the global well-posedness using the first non trivial invariant) in $H^1(\R).$ The new ingredient is to apply a gauge transformation (a variant of the classical Cole-Hopf transform for the Burgers equation) in order to eliminate the terms involving the interaction of very low and very high frequencies, where the derivative falls on the very high frequencies. Note that those interactions are responsible of the lack of regularity of the flow map, see last paragraph and \cite{KT2}.  

More precisely, to obtain the estimates at a $H^1$ level, Tao introduces the new unknown

$$w=\partial_xP_{+\text{hi}}(e^{-1/2F}),$$

where F is some spacial primitive of u and $P_{+\text{hi}}$ is the projection on high positive frequencies. Then w satisfies the equation

$$\partial_tw-i\partial_x^2w=-\partial_xP_{+\text{hi}}(\partial_x^{-1}wP_-\partial_xu)+\text{negligible terms},$$

where $P_-$ is the projection on negative frequencies.

Thanks to the frequency projections, the nonlinear term appearing in the right hand side does not involve any low-high frequency interaction terms. Finally, to invert this gauge transform, one gets an equation of the form 

\begin{equation}\label{tao}
u=2ie^{\frac{i}{2}F}w+\text{negligible terms}.
\end{equation}

The gauge transformation is not related to the complete integrability of the BO equation since a variant of the gauge transform is used in \cite{HIKK} in order to establish the global well-posedness in $L^2$ of the (non integrable) fractional KdV equation

\begin{equation}\label{fKdV}
u_t+uu_x-D^\alpha u_x=0
\end{equation}

for $\alpha\in (1,2).$

Further improvements were brought by Burq and Planchon \cite{BP} who proved LWP in $H^s(\R), s>1/4$ and thus global-wellposedness in the energy space $H^{1/2}(\R).$

Finally Kenig and Ionescu \cite{IK} proved the local (and thus global) well-posedness in the space $H^s(\R), s\geq 0.$ Uniqueness however is obtained in the class of limits of smooth functions. Both \cite{IK} and \cite{BP} use Tao's ideas in the context of Bourgain spaces. The main difficulty is that Bourgain's spaces do not enjoy an algebra property so that one looses regularity when estimating u in terms of w in \eqref{tao}. To overcome this difficulty Burq  and Planchon first paralinearize the equation and then use a localized version of the gauge transform on the worst nonlinear term. On the other hand, Ionescu and Kenig decompose the solution in two parts : the first one is the smooth solution of BO evolving from the low frequency part of the initial data while the second one solves a dispersive system renormalized by a gauge transformation involving the first part. This system is solved via a fixed point argument in a dyadic version of Bourgain's space with  a special structure in low frequencies.

In \cite{MP} Molinet and Pilod simplified the proof of Ionescu and Kenig and furthermore proved stronger uniqueness properties, in particular they obtain the unconditional uniqueness in the space $L^\infty(0,T; H^s(\R))$ for $s>1/4.$

Methods using gauge transformations are very good to obtain low regularity results but behave badly with respect to perturbations of the BO equation. In particular it is not clear if they apply to the ILW equation. 

Molinet and Vento \cite{MV} proposed a method which is less powerful to get low regularity results (say in $L^2(\R)$) but allows to deal with perturbations of the BO equation, in particular to the ILW equation. Their approach combines classical energy estimates with Bourgain type estimates on a time interval that does not depend on the space frequency. It yields local well-posedness in $H^s(\R), \geq 1/2,$ the unconditional uniqueness holding only when $s>1/2.$ The method  works as well in the periodic case and also for more general dispersion symbols, in particular similar  results apply to the ILW equation. Since it does not rely on the gauge transform the method allows to prove strong convergence results in the energy space for solutions of viscous perturbations of the BO equation (see {\it eg} \cite {OS} for physical examples and \cite {Mat13} for a formal multi-soliton dissipative perturbation of BO).

 Combining this technique with refined Strichartz estimates and modified energies, Molinet,  Pilod and Vento \cite{MPV} extended this result to the regularity $H^s(\R), s>\frac{1}{4}.$ Moreover their method also applies to fractional KdV equation with low dispersion,  that is to 

$$u_t+uu_x-D^\alpha u_x=0,$$

with $0<\alpha\leq 1$ yielding a local well-posedness result in $H^s(\R), s>\frac{3}{2}-\frac{5\alpha}{4},$ and thus the global well-posedness in the energy space $H^{\alpha/2}(\R)$ for $\alpha>\frac{6}{7}.$

Recently Ifrim and Tataru \cite{IT} revisited the global well-posedness of BO in $L^2(\R)$ by  a normal form approach. More precisely they split the quadratic part into two parts, a milder one and a paradifferential one. The normal form correction is then constructed in two steps, a direct quadratic correction for the milder part and   a renormalization type correction for the paradifferential part. The second step use a paradifferential version of Tao's renormalization. 

\subsection{Long time dynamics}
We have seen that the Cauchy problem is well understood in Sobolev spaces $H^s(\R), s\geq 0.$ The long time dynamics is much less understood (see subsection 5.4 below for conjectures on solutions emerging from initial data of arbitrary size). The case of small initial  data is considered in \cite{IT} where the following result is proved

\begin{theorem}\label{IfTa}
Assume that the initial data $u_0$ for the BO equation satisfies

$$|u_0|_2+|xu_0|_2\leq \epsilon \ll 1.$$

Then the global associated solution u satisfies the dispersive decay bounds

$$|u(t,x)|+|Hu(x,t)|\lesssim \epsilon|t|^{-1/2}\langle x_-t^{-1/2}\rangle^{-1/2}, \;\text{where}\; x_-=-\min (x,0)$$
up to time

$$|t|\lesssim T_\epsilon:=e^{\frac{c}{\epsilon}},\quad c\ll 1.$$
\end{theorem}
\begin{remark} We recall that  the linear BO flow satisfies the decay property

$$||e^{-tH\partial_x^2}||_{L^1\to L^\infty}\lesssim t^{-1/2}.$$

The better decay rate in the region $x<0$ is due to the positivity of the phase velocity $|\xi|$ which send the propagating waves to the right and the dispersive waves to the left.
\end{remark}

We refer to \cite{IS} for similar results on ILW.

\subsection{Solitary waves}
Both the BO and ILW equations possess explicit solitary wave solutions. As aforementioned, Benjamin \cite{Ben} found the profile of the BO one, namely
\begin{equation}\label{SWBO}
\phi(x)=\frac{4}{1+x^2},
\end{equation}

leading to the family of solitary waves $\phi_c(x-ct)$ where c is a positive constant and 

$$\phi_c(y)=\frac{4c}{c^2y^2+1}.$$

Note that since the symbol of the dispersion operator is not smooth, Paley-Wiener type arguments imply that the solitary wave cannot be exponentially decreasing.

Note also that $|\phi_c|_1\equiv 4\pi, \forall c>0$ while $|\phi_c|_2, |\phi_c|_\infty \to 0$ as $c\to 0,$ thus BO possesses arbitrary small solitary waves.

On the other hand the symbol in the ILW equation is smooth and the explicit solitary wave found in \cite{J} decays exponentially, namely for arbitrary $C>0$ and $\delta>0$ 

\begin{equation}\label{SWILW}
\phi_{C,\delta}(x)=\frac{2a\sin (a\delta)}{\cosh (ax)+\cos(a\delta)},
\end{equation}

where $a$ is the unique solution of the transcendental equation 

$$a\delta\cot(a\delta)=(1-C\delta),\quad a\in(0,\pi/\delta).$$

Before considering the stability properties of those solitary waves, we review important {\it uniqueness} results. Amick and Toland (\cite{AmTo2}, see also \cite{AmTo, Alb1}), using the maximum principle for linear elliptic equations, estimates on a Green function and the Cauchy-Riemann equation, proved the uniqueness of the BO solitary wave in the sense that the only solutions to the pseudo-differential equation

$$u(x)^2-u(x)=\mathcal G(u)(x), \quad x\in \R,$$
which satisfies the boundary condition 

$$u(x)\to 0\; \text{as}\; |x|\to \infty,$$

where
$$\mathcal G(f)(x)=\frac{1}{2\pi}\int_{-\infty} ^\infty|\xi|e^{-i\xi x}\left(\int_{-\infty}^\infty f(\eta)e^{i\xi \eta}d\eta\right) d\xi$$

are the functions

$$u(x)=0\quad \text{and}\quad u_a(x)=\frac{4}{1+(x-a)^2}, \; a\in \R.$$

\vspace{0.3cm}
A similar result holds true for the ILW equation. Let 

\begin{equation}\label{eqSWILW}
(\mathcal N_\delta+\gamma)\phi= \phi^2
\end{equation}

where \footnote{Note that $\mathcal F(\mathcal N_\infty u)(\xi)=|\xi|\hat u (\xi).$}

$$\mathcal F(\mathcal N_\delta u)(\xi)=(\xi\coth\xi\delta)\mathcal F(u)\xi$$ 

be the equation of a solitary wave $ u(x,t)=\phi(x-Ct)$ to the ILW equation, where $\gamma= C-1/\delta.$

Then Albert and Toland (\cite{AlTo} and also \cite{Alb1}) proved that for $\delta>0$ and $C>0$ be given, if $\phi\in L^2(\R)$ is a non trivial solution of \eqref{eqSWILW}, then  there exists $b\in \R$ such that $\phi(x)=\phi_{C,\delta}(x+b).$

The proof in \cite{Alb1} relies on two special properties of the operator $\mathcal N _\delta$ that degenerate into classical ones (linked to the Hilbert transform) when $\delta =\infty,$ providing another proof of the Amick-Toland uniqueness result for the BO equation.

\vspace{0.3cm}
We turn now to stability issues.  The orbital stability of the BO and ILW solitary waves can be obtained by the classical Cazenave-Lions method in \cite{CZ} (minimization of the Hamiltonian with fixed $L^2$ norm) providing orbital stability in the energy space $H^{1/2}(\R)$ but the first known proofs were by using the Souganidis-Strauss method, see \cite {BBSSB, ABH, BoSo, BSS} for BO and \cite{AB1, ABH, Alb2} for ILW.

The $H^1$ orbital stability of the Benjamin-Ono 2-soliton is proven in \cite{NL}.
The 2-soliton of velocities $c_1>0, c_2>0$ with $c_1<c_2$ is explicitly given in \cite{Mat2}) as 

$$\phi_{c_1,c_2}(t,x)=4\frac{c_2\theta_1^2+c_1\theta_2^2+(c_1+c_2)c_{12}}{(\theta_1\theta_2-c_{12})^2+(\theta_1+\theta_2)^2},$$

where $\theta_n=c_n(x-c_nt), n=1,2$ and $c_{12}=\left(\frac{c_1+c_2}{c_1-c_2}\right)^2.$

The proof in \cite{NL} relies on the integrability of BO, and is reminiscent of the similar one for the stability of the KdV 2-soliton (see \cite{MaSa}), namely it uses the fact that the 2-soliton locally minimize the invariant $I_4$ subject to given values of $I_3$ and $I_2.$

An alternative characterization of the 2-soliton uses the self-adjoint operator M that appears in the Lax pair associated to the integrable equation (KdV or BO). Recall that the Lax pair for BO was constructed in \cite{AF}. The 2-solitons are the potentials for which the self-adjoint operator M has two eigenvalues.

The proof involves the spectral analysis of the one-parameter family of self-adjoint operators $L(t)$ which are the linearization of 

$$I_4'(u)+\alpha I_3'(u)+\beta I_2'(u)=0$$

at the double soliton.
Contrary to the KdV case where $L(t)$ is a fourth order self-adjoint linear ordinary differential operator, for BO one obtains a nonlocal operator since the Hilbert transform appears in it. This makes the spectral analysis more complicated. The new approach in \cite{NL} consists in making a simplification in the spectral problem to reduce the spectral analysis of the one-parameter family $L(t)$ to the analysis of the spectra of two stationary operators $L_1$ and $L_2$. The proof is then reduced to proving the two facts:

${\bf 1.}\; L_1$ has one negative eigenvalue and $L_2$ has no negative eigenvalue;

${\bf 2.}$ zero is a simple eigenvalue of $L_1$ and $L_2.$

\vspace{0.3cm}
The asymptotic stability of the BO solitons in the energy space $H^{1/2}(\R)$ is proven in  \cite{GTT, KeMa}.

We describe briefly the statement of \cite{KeMa}, denoting by 

$$Q(x)=\frac{4}{1+x^2}$$

the profile of the BO soliton and $Q(cx)= cQ(cx). $
\begin{theorem}\label{multBO}
There exist $C,\alpha_0>0,$ such that if $u_0\in H^{\frac{1}{2}}(\R)$ satisfies $||u_0-Q_c||_{\frac{1}{2}}=\alpha \leq\alpha_0,$ then there exist $c^+>0$ with $|c^+-1|\leq C\alpha$ and a $C^1$ function $\rho(t)$ such that the solution of BO with $u(0)=u_0$ satisfies

$$u(t,.+\rho'(t))\to Q_{c^+}\quad \text{in}\;H^{\frac{1}{2}}\; \text {weak},\quad |u(t)-Q_{c^+}(.-\rho(t))|_{L^2(x>\frac{t}{10})}\to 0,$$

$$\rho(t)\to c^+\quad \text{as}\; t\to +\infty.$$
\end{theorem}

\begin{remark}
The convergence of $u(t)$ to $Q_+$ as $t\to \infty$ holds in fact strongly in $L^2$ in the region $x>\epsilon t,$ for any $\epsilon >0$ provided $\alpha_0=\alpha_0(\epsilon)$ is small enough. This result is optimal in $L^2$ since $u(t)$ could contain other small (and then slow) solitons and since in general $u(t)$  does not go to $0$ in $L^2$ for $x<0.$ 

For instance, if $||u(t)-Q_{c^+}(.-\rho(t))||_{\frac{1}{2}}\to 0$ as $t\to \infty$, then $E(u)=E(Q_{c^+})$ and moreover $\int_\R u^2dx=\int_\R Q_{c^+}^2dx$ so that by the  variational characterization of $Q(x)$, $u(t)=Q_{c^+}(x-x_0-c^+t)$ is a solution.

It is expected but not proved yet that the convergence  in the same local sense $x>\epsilon t$ holds in $H^{\frac{1}{2}}(\R).$
\end{remark}

The proof of Theorem 5 is based on the corresponding one  (\cite{MaMe} and the references therein) for the generalized KdV equation
 where the stability is deduced from a Liouville type theorem. There are however two new difficulties.

1. The proof of the $L^2$ monotonicity property is more subtle because of the nonlocal character of the BO equation. 

 2. The proof of the linear Liouville theorem which requires the analysis of some linear operators related to $Q.$

Similar arguments and the strategy used for the asymptotic stability in the energy space of the sum of N solitons for the subcritical generalized KdV equation yield a similar result for N-solitons of the BO equation (\cite{MaMe}):

\begin{theorem}
Let $N\geq 1$ and $0<c_1^0<...<c_N^0.$ There exist $L_0>0, A_0>0$ and $\alpha_0>0$ so that if $u_0\in H^{1/2}$ satisfies for some $0\leq\alpha<\alpha_0, L\geq L_0,$

$$||u_0-\sum_{j=1}^NQ_{c_j^0}(\cdot-y_j^0)||_{H^{1/2}}\leq \alpha \quad\text{where}\quad \forall j\in (2,...N),\;y_j^0-y^0_{j-1}\geq L,$$

and if $u(t)$ is the solution of BO corresponding to $u(0)=u_0,$ then there exist $\rho_1(t),...,\rho_N(t)$ such that the following hold 

(a) Stability of the sum of N decoupled solitons,

$$\forall t\geq 0,\; ||u(t)-\sum_{j=1}^NQ_{c_j^0}(x-\rho_j(t)||^{1/2}\leq A_0\left(\alpha+\frac{1}{L}\right).$$

(b)Asymptotic stability of the sum of N solitons. There exist $c_1^+,...,c_N^+,$ with $|c_j^+-c_j^0|\leq A_0(\alpha+\frac{1}{L})$ such that

$$\forall j,\quad u(t,.+\rho_j(t)))\rightharpoonup Q_{c_j^+} \quad \text{in}\; H^{1/2}\;\text{weak as}\; t\to +\infty,$$

$$||u(t)-\sum_{j=1}^NQ_{c_j^+}(.-\rho_j(t)||_{L^2(x\geq\frac{1}{10c^0_1(t)})}\to 0, \; \rho'_j(t)\to c_j^+ \;\text{as}\; t\to +\infty.$$

\end{theorem}

Recall that the BO equation possesses explicit multi-soliton solutions. Let $U_N(x;c_j, y_j)$ denotes the explicit family of N-soliton profiles. Using the previous theorem and the continuous dependence of the solution in $H^{1/2}$ one obtains the following corollary \cite{MaMe}.

Let $N\geq 1, 0<c_1^0<...<c_N^0$ and set

$$d_N(u)=\inf\lbrace||u-U_N(.;c_j^0,y_j)||_{1/2}, y_j\in \R \rbrace.$$

Then 

\begin{corollary} (Stability in $H^{1/2}$ of multi-solitons).
For all $\delta >0,$ there exists $\alpha >0$ such that if $d_N(u_0)\leq \alpha$ then for all $t\in \R,$ $d_N(u(t))\leq \delta.$
\end{corollary}

As aforementioned, the proof of Theorem \ref{multBO} is based on a rigidity theorem:

\begin{theorem}
There exist $C,\alpha_0>0,$ such that if $u_0\in H^{\frac{1}{2}}$ satisfies $||u_0-Q_c||_{H^{\frac{1}2}}=\alpha \leq\alpha_0,$ and if the solution $u(t)$ of BO with $u(0)=u_0$ satisfies for some function $\rho(t)$

$$\forall\epsilon>0, \exists A_\epsilon>0, \;\text{such that}\; \int_{|x|>A_\epsilon} u^2(t,x+\rho(t))dx<\epsilon,$$
 then there exist $c_1>0, x_1\in \R$ such that

 $$u(t,x)=Q_{c-1}(x-x_1-c_1t),\quad |c_1-1|+|x_1|\leq C\alpha.$$

\end{theorem}

The proof of the rigidity theorem (Liouville theorem) requires in particular the analysis of some linear operators related to $Q$ that uses the fact that $Q(x)$ is explicit and some known results on the linearization of the BO equation around $Q$ \cite{BBSSB, We}.

\begin{remark}
Similar asymptotic stability results - though expected- do not seem to be known for the ILW equation.
\end{remark}

\subsection{A result on long time asymptotic}

As already noticed, the complete description asymptotic behavior of the BO or ILW solutions with arbitrary large initial data is unknown (see Section 5.4 for the so-called soliton resolution conjecture).

A significant progress was made in \cite{MuPo} which is concerned with the long time behavior of solutions. It concerns global solutions u of  the Benjamin-Ono equation  satisfying

\begin{equation}\label{MP0}
 u\in C(\R;H^1(\R))\cap L^\infty_{\text{loc}}(\R; L^1(\R)),
\end{equation}
and moreover :

$$\exists a\in [0,1/2),\; \exists c_0>0\quad \text{such that}$$

\begin{equation}\label{hypCl}
 \forall T>0, \;\sup_{t\in[0,T]}\int_{-\infty}^\infty|u(x,t)|dx\leq c_0(1+T^2)^{a/2}.
\end{equation}

The main result in \cite{MuPo} is then

\begin{theorem}\label{MP}
Under the above assumption,
\begin{equation}\label{MP1}
\int_{10}^\infty \frac{1}{\log t}\left(\int_{-\infty}^\infty \phi'\left(\frac{x}{\lambda(t)}\right)\left(u^2+(D^{1/2}u)^2\right)(x,t)dx\right)dt<\infty.
\end{equation}
Hence,
\begin{equation}\label{MP2}
\liminf_{t\to +\infty}\int_{-\infty}^\infty  \phi'\left(\frac{x}{\lambda(t)}\right)\left(u^2+(D^{1/2}u)^2\right)(x,t)dx=0,
\end{equation}

with 
\begin{equation}\label{MP3}
\lambda(t)=\frac{ct^b}{\log t},\quad a+b=1,\quad \text{and}\; \phi'(x)=\frac{1}{1+x^2},
\end{equation}
for any fixed $c>0.$
\end{theorem}

\begin{remark}
1. This result discards the existence of non trivial time periodic solutions (in particular breathers) and of solutions moving with a speed slower that a soliton.

2. The theorem implies that there exists a sequence of times $\lbrace t_n; n\in \mathbb{N}\rbrace$ with $t_n\to +\infty$ as  $n\to \infty$ such that
\begin{equation}\label{MP3}
\lim_{n\infty}\int_{|x|\leq \frac{ct_n^b}{\log (t_n)}} (u^2+(D^{1/2}u)^2)(x,t_n)dx=0.
\end{equation}

3. The solitons 

$$u(x,t)=Q_c(x-ct), \quad Q_c(x)=\frac{4c}{1+c^2x^2},$$
belongs to the class \eqref{MP0} and they also satisfy \eqref{MP2}.
\end{remark}

\begin{remark}
The above result does not use the integrability of the BO equation. It is likely to hold for the ILW equation.
\end{remark}

The proof of Theorem \ref{MP} relies in particular on the estimates obtained in \cite{KeMa} and uses the fact that $\phi'$ is a multiple of the soliton.
\section{Rigorous results by IST methods}

\subsection{The BO equation}
An attempt to solve the Cauchy problem of the BO equation by IST methods is due to Anderson and Taflis \cite{AT}. They obtained a formal series linearization of the BO equation that can be interpreted as a distorted Fourier transform associated to a singular perturbation of $\frac{1}{i}\frac{d}{dx},$ from which they deduce  a Lax pair for BO and under suitable assumptions on the scattering data, a power series for the inverse transform. Unfortunately, as proved in \cite{CW}, both the direct and inverse  problems considered in \cite{AT} are not analytic for generic potential, leading to divergent series in general.

The first rigorous results on the Cauchy problem by IST methods were obtained by Coifman and Wickerhauser in \cite{CW}, (see also  \cite{BC1, BC2}). They used a more complicated regularized IST formalism and solved it by iteration, yielding the global well-posedness of the Cauchy problem for small initial data. More precisely they use constructive methods to investigate the spectral theory of the Benjamin-Ono equation. Since the linearization series used previously is singular (see above), they replace it with an improved series obtained by finite-rank renormalization. This introduces additional scattering data, which are shown to be dependent upon a single function, though not the usual one. They then prove the continuity of the direct and inverse scattering transforms defined by the improved series for small complex potentials. 

Let $w(x)=(1+|x|).$ The global well-posedness in \cite{CW} writes :
\begin{theorem}\label{CW}
Let $q$ be a real function such that $w^{n+1}q$ is a small function in $L^1(\R)$ for some $n>0.$ Suppose also that $w^nq'$ and $w^nq"$ are also small in $L^1(\R).$ Then there exists a unique solution of the Benjamin-Ono equation with initial data $q.$
\end{theorem}

Some interesting issues arise also in \cite{CW} linked to the generation of soliton solutions. We comment them briefly in connection with \cite{ Mat9,MPST} and specially \cite{ PS}.
Let the initial data for the BO equation renormalized as  $u_0=U_0U(x/L_0),$ where $U_0$ is a characteristic amplitude and $L_0$ is a characteristic length that we will follow closely. Then, the number of solitons depends on the sole parameter $\sigma=U_0L_0$ (Ursell number). In the limit $\sigma \gg1$ the initial potential generates a large number $N$ of solitons. An approximation for $N$ was found in \cite{Mat9} and later confirmed by \cite{MPST},

$$N=\frac{1}{2\pi}\int_{u_0(x)\geq 0} u_0(x)dx.$$
An important question is that of the existence of a threshold for the generation of solitons, that is whether a small initial perturbation of an  algebraic profile  $u(x)=\frac{a}{1+x^2}$ {\it ie} in the limit 
$\sigma \ll1, (a\ll 1)$ can support propagation of at least one soliton.

Note that the mass $\mathcal M [u_c]=\int_\R u_c(x)dx$ of the BO soliton $u_c(x-ct)=\frac{4c}{1+c^2(x-ct)^2}$ is constant,  $\mathcal M [u_c]=4\pi.$

For the modified KdV equation, such a property is related to the existence of a threshold on the soliton generation, in this case perturbations with $\mathcal M[u]\leq \frac{1}{2}\mathcal M _{\text{sol}}$ do not support solitons.

For the BO equation, this issue depends on the possible non genericity of the potentials and on the very special structure of the Jost function.  By definition generic potentials $u$ are those for which  $n_0\neq 0$ where

$$n_0=\frac{1}{2\pi}\int_\R u(x)n(x)dx.$$

Here $n(x)$ is the limiting Jost function that we describe now. For generic potentials, the Jost functions vanish as $k\to 0^+$ according to the approximation

$$N(x,k)\to \frac{n(x)}{1+n_0(\gamma+\ln(ik))}+O\left(\frac{k}{\ln k}\right),$$

where $\gamma$ is the Euler constant and $n(x)$ satisfies

$$in_x=-P^+(un_x)$$

where $P^+(v)=\frac{1}{2}(v-iHv).$

The problem at the center of the analysis is the spectral problem

\begin{equation}\label{specBO}
i\phi_x^++k(\phi^+-\phi^-)=-u(x)\phi^+.
\end{equation}

Coifman and Wickerhauser \cite{CW} proved that the scattering problem \eqref{specBO} has no bound states in a  neighborhood of the origin, if $u(x)\to O(|x|^{-1-\mu}), \mu>0.$ In fact this result is valid for generic potentials in the previous sense.

For non generic potentials, {\it ie} satisfying $n_0=0,$ the limiting Jost function is bounded in the limit $k\to 0^+$
 and properties of the scattering problem are modified. The zero potential $u(x)=0$ is non generic as well as the soliton solutions.

 The main results in \cite{PS} concerns the study of the perturbation of non generic potentials where the number of bound states may change depending on the type of the perturbation. More precisely the authors consider a potential in the form $u^\epsilon=u(x)+\epsilon \eta(x)$ where $\epsilon\ll 1$ and $u(x)$ satisfies the constraint $n_0=0.$ in the particular case $u(x)=0,$ this reduces to the problem of soliton generation by a small initial data.  They derive a criterion for a new eigenvalue to emerge from the edge of the continuum spectrum at $k=0.$ In particular, new eigenvalues may always appear due to perturbations of the zero background and the soliton solutions.  They derive the leading order term  for the new eigenvalue and for the associated bound state at short and long distances, as well as for the variation of the continuous spectrum.

 \vspace{0.3cm}
Further progress of the IST theory for BO is due to Miller and Wetzel \cite{MW2} who studied the direct scattering problem of the Fokas-Ablowitz IST theory when the potential is a rational function with simple poles and obtained the explicit formula for the scattering data.

\vspace{0.3cm}
A  breakthrough on the IST for the BO equation appeared recently in the works by Yulin Wu \cite{Wu, Wu2} who in particular solved completely the direct scattering problem for data of arbitrary size.

In \cite{Wu} Wu studies the operator $L_u$ arising in the Lax pair for the BO equation and proved that its spectrum is discrete and simple. More precisely, we first rewrite the Fokas-Ablowitz formalism in a slightly different way, denoting first

$$C_\pm \phi=\frac{\phi\pm iH\phi}{2}$$

the Cauchy projection.

 When they act on $L^2(\R)$ the ranges are $H^\pm$, the Hardy spaces of $L^2$ functions whose Fourier transform are supported on the positive and negative half lines. The Lax pair writes in those notations, on $H^+$ 
 
 $$L_u\phi=\frac{1}{i}\phi_x-C_+(uC_+\phi)=\lambda \phi$$
 
 $$B_u\phi =2\lambda \phi_x+i\phi_{xx}+2(C_+u_x)(C_+\phi),$$
 
 and in $H^-:$
 
 $$L_u\phi=\frac{1}{i}\phi_x-C_-(uC_-\phi)=\lambda \phi$$
 
 $$B_u\phi =2\lambda \phi_x+i\phi_{xx}+2(C_-u_x)(C_-\phi).$$
 
Since when $u$ is real the equations on $H^-$ are just the complex conjugate of the equations on $H^+$ one will assume that $u$ is real and focus on the $H^+$ part of the Lax pair. The scattering data of the IST are closely related to the spectrum of the operator $L_u.$

The next two theorems in \cite{Wu} provide useful informations on the spectral properties of $L_u.$

\begin{theorem}\label{Wu1}
Suppose that $u\in L^2(\R)\cap L^\infty(\R).$ Then $L_u$ is  a relatively compact perturbation of $\frac{1}{i}$ and is self-adjoint on $H^+$ with domain $H^+\cap H^1(\R).$
\end{theorem}

\begin{theorem}\label{Wu2}
Suppose that $u\in L^1(\R)\cap L^\infty(\R)$ and $xu\in L^2(\R).$ Then the operator $L_u$ has only finitely many negative eigenvalues and the dimension of each eigenspace is $1.$ 
\end{theorem}

\begin{remark}
By Weyl's theorem and standard spectral theory, the essential spectrum of $L_u$ is the same as that of $\frac{1}{i}\partial_x$, that is $\R^+\cup \lbrace 0 \rbrace.$
However it is not clear that the eigenvalues are simple and even if there are finitely many of them. On the other hand those spectral properties are crucial for the construction of scattering data in the Fokas-Ablowitz  IST method in \cite{AF}.

A crucial step in the  proof of simplicity is the discovery of  a new identity connecting the $L^2$ norm of the eigenvector to its inner product with the scattering potential. 

More precisely, assuming that $u\in L^2(\R)\cap L^\infty(\R),$ and $xu\in L^2(\R),$, then if $\lambda <0$ is an eigenvalue of $L_u$ and $\phi$ is an eigenvector, then

$$\left|\int_\R\phi dx\right|^2=2\pi\lambda\int_\R |\phi|^2ds.$$

The proof for finiteness is an extension of ideas involved in the Bergman-Schwinger bound for Schr\"{o}dinger operators.
\end{remark}

\vspace{0.3cm}
As aforementioned, the direct scattering problem is completely justified in \cite{Wu2} where Alen Wu examined the full spectrum of $L_u$ establishing existence, uniqueness and asymptotic properties of the Jost solutions to the scattering problem, providing possible directions for the correct setup for the inverse problem.

 Recalling that 
$w(x)=(1+|x|),$ we will use the notations

 $L^p_s(\R)=\lbrace f, w^s f\in L^p(\R)\rbrace$ and 

$$H_s^s(\R)=\lbrace f; f\in L^2_s(\R) \text{and} \;\hat f \in L^2_s(\R)\rbrace, $$

Two Jost functions $m_1$ and $m_e$ are considered in \cite{Wu2}. They are solutions of the following equations with suitable boundary conditions stated in the next lemma :

$$\frac{1}{i}\partial _xm_1-C_+(um_1)=k(m_1-1),$$

and 

$$\frac{1}{i}\partial _xm_e-C_+(um_e)=\lambda m_e.$$  

Here  $\lambda\pm 0i\in \R^+\pm 0i,$  and

$$k\in\rho(L_u)\cup(\R^+\pm 0i)=(\mathbb{C}\setminus \lbrace \lambda_1....\lambda _N\rbrace\setminus [0,\infty))\cup(\R^+\pm 0i),$$

which is the resolvent set glued with two copies of the positive real axis. Following again \cite{Wu2} we provide here the translation of   notation:

$$M(x,k)=m_1(x, \lambda+ 0i), \bar{M}(x,\lambda)=m_e(x, \lambda +0i),$$

$$N(x,\lambda)=m_e(x,\lambda-0i),\bar{N}(x,\lambda)=m_1(x,\lambda-0i).$$

The following lemma may be considered as the definition of $m_1$ and $m_e.$ In what follows, one uses the integral operators

$$G_k(x)=\frac{1}{2\pi}\int_0^\infty \frac{e^{ix\xi}}{\xi-k} d\xi,$$

for $k\in\mathbb{C}\setminus [0,\infty),$ and

$$\tilde{G}_k(x)=\frac{1}{2\pi}\int_{-\infty}^0 \frac{e^{ix\xi}}{\xi-k} d\xi,$$


for $k\in \mathbb{C}\setminus (-\infty,0].$

 For $\epsilon >0$ one has 
 
 \begin{equation}\label{38}
 G_{\lambda\pm 0i}(x)=\frac{1}{2\pi}\int_{-\infty}^\infty \frac{e^{is\xi}}{\xi-(\lambda\pm i\epsilon)}d\xi-\tilde{G}_{\lambda\pm i\epsilon}(x)=\pm ie^{\mp \epsilon x}e^{i\lambda x}\chi_{\R^{\pm}}(x)-\tilde{G}{\lambda \pm i\epsilon}(x),
 \end{equation}
 
 with 
 
 \begin{equation}\label{39}
 G_{\lambda\pm 0i}(x)=\lim_{\epsilon \to 0}G_{\lambda\pm i\epsilon}(x)=\pm ie^{i\lambda x}\chi_{\R^{\pm}}(x)-\tilde{G}\lambda (x),
 \end{equation}
 
 for $\lambda >0.$
 
 The limit in \eqref{39} holds in the following sense  : the first term in \eqref{38} converges pointwise and the second term in \eqref{38} converges in  $L^{p'}$ for every $p'\in [2,\infty).$ The latter is checked when observing that $\tilde{G}_{\lambda\pm i\epsilon}$ is the inverse Fourier transform of
 
 $\frac{\chi_\R-\xi}{\xi-(\lambda\pm i\epsilon)},$ which converges to $\frac{\chi_\R-\xi}{\xi-\lambda}$ in every $L^p$ for  $p\in (1,2]$ assuming  $\lambda >0.$

\begin{lemma}
Let $p>1$ and $s>s_1>1-\frac{1}{p}$ be given and let $u\in L^p_s(\R).$ Assume that $m_1(x,k),m_e(x,\lambda \pm 0i)\in L^\infty_{-(s-s_1)}(\R)$ for fixed $k\in(\mathbb{C}\setminus[0,\infty))\cup(\R^+\pm 0i)$ and $\lambda\in \R^+,$  then the following are equivalent:

(a) $m_1(x,k), m_e(x,\lambda\pm 0i)$ solve

$$\frac{1}{i}\partial_x m_1-C_+(um_1)=k(m_1-1),$$

$$\frac{1}{i}\partial_x m_e-C_+(um_e)=\lambda m_e,$$

together with the asymptotic conditions
\begin{equation}
m_1(x,k)-1\to 0\; \begin{cases}
\text{as}\;  |k| \to \infty &\text{if}\; k\in \mathbb{C}\setminus [0,\infty)\\
\text{as} \; x\to \mp \infty &\text{if}\; k=\lambda \pm 0 i\in \R^+\pm 0 i,
\end{cases}
\end{equation}
$$m_e(x,\lambda\pm 0i)-e^{i\lambda x}\to 0 \quad \text{as}\; x\to \mp \infty.$$

The above asymptotic conditions should be read with either the upper or the lower sign.

(b) $m_1(x,k), m_e(x,\lambda \pm 0i)$ solve the following integral equations :

\begin{equation}\label{314}
m_1(x,k)=1+G_k \star (um_1(.,k))(x),
\end{equation}

\begin{equation}\label{315}
m_e(x,\lambda\pm 0i)=e^{i\lambda x}+G_{\lambda\pm 0i)}\star(um_e(.,\lambda\pm 0i))(x).
\end{equation}

Moreover, if either (a) or (b) holds, one has the stronger bounds

$$m_1(x,k)-1\in L^\infty(\R)\cap \mathbb{H}^{p,+}$$

for fixed $k\in \mathbb{C}\setminus [0,\infty)$ and 

$$m_1(x,\lambda\pm 0i),m_e(x,\pm 0i)\in L^\infty(\R)$$

for fixed $\lambda \in \R^+.$
\end{lemma}

Here is the theorem establishing the existence and uniqueness of Jost solutions, as stated in \cite{Wu2}.

\begin{theorem}\label{Jost}
Let $s>s_1>\frac{1}{2}$ and $u\in L^2_s(\R).$ Let $\rho(L_u)$ be the resolvent set of  $L_u=\frac{1}{i}\partial_x-C_+uC_+,$ regarded as an operator
 on $\mathbb{H}^+.$
 Then for every $k\in \rho(L_u)\cup (\R^+\pm 0i),$ and every $\lambda >0,$ there exist unique $m_1(x,k)$ and $m_e(x,\lambda\pm 0i)\in L^\infty_{-(s-s_1)}(\R)$ solving   
 \eqref{314}, \eqref{315}
 respectively, with improved bounds $m_1(x,k), m_e(x, \lambda\pm 0i)\in L^\infty(\R).$  Furthermore, the mapping $k\to m_1(k)$ is analytic from $\rho(L_u)$ to $L^\infty_{-(s-s_1)},$ and $m_1(k)\in C^{0,\gamma}_{\text{loc}}((\rho(L_u)\cup \R^+\pm 0i)),L^\infty_{-(s-s_1)}(\R)$ while $m_e(\lambda\pm 0i)\in C^{0,\gamma}_{\text{loc}}(\R).$

Here $\gamma$ is some number between $0$ and $1.$

\end{theorem}

When one  studies the asymptotic behavior of the Jost solutions and scattering coefficients as k approaches  to $0$ within the set $\rho(L_u)\cup(\R^+\pm0i)$ one notices that the convolution kernel $G_k$ has a logarithmic singularity at $k=0$ and so does the operator $T_k.$ By substracting a rank one operator from $T_k$  the modified operator has  a limit at $k=0.$ The asymptotic behavior of the Jost functions can be recovered from  {\it modified} Jost functions, we refer to \cite{Wu2} for details. One recovers rigorously the asymptotics in \cite{AF,KaMat}. Those asymptotic formulas are useful to clarify the global behavior of the scattering coefficients. One also finds in \cite {Wu2} asymptotic formulas in the limit $k\to \infty$ and a formal derivation of the time evolution of the Jost functions and scattering coefficients provided $u$ is a smooth and decaying solution of  the BO equation.





\vspace{0.3cm}
\begin{remark} A large class of {\it quasi-periodic} solutions of the BO equation has been found in \cite{SI} by the Hirota method and in \cite{DK} by the method of the theory of finite-zone integration.  More precisely the quasi-periodic solution of the BO equation found in \cite{DK} writes 

\begin{equation}\label{DK}
u(x,t)=c+\sum_{j=1}^N (a_j-b_j)-2\; \text{Im}\; \frac{\partial}{\partial x} \ln \det M(x,t),
\end{equation}

where the $N\times N$ matrix $M(x,t)$ has the elements

\begin{equation}\label{matrix}
M_{jm}=c_m\delta_{jm}e^{i(a_m-b_m) -i(a_m^2-b_m^2)t}-\frac{1}{b_j-a_m}.
\end{equation}

The constants $a_j, b_j, c$ satisfy the inequalities $c\leq a_1\leq b_1\leq a_2\leq b_2\leq ...\leq a_n\leq b_n$ and the constants $c_i$ are defined by

\begin{equation}\label{defc}
|c_i|^2=-(b_i-c_i)\prod_{j\neq i}^N\frac{(a_i-a_j)(b_i-b_j)}{(a_i-c)\prod_{j=1}^N(b_i-a_j)(a_i-b_j)}.
\end{equation}
This multiply periodic N-phase solution has been used in \cite{JMS} to construct the modulation solution of BO corresponding to  a step like function:
\begin{equation}
u(x,0)= \begin{cases} A\quad x<0,\;A>0\\
0 \quad x>0,
\end{cases}
\end{equation}

or a smoothly decreasing  function (leading to a shock formation for the Burgers equation):

\begin{equation}
u(x,0)=\alpha \left(1-\frac{2}{\pi}\arctan \beta x\right),\quad \alpha, \beta >0,
\end{equation}

or a modulated wavetrain which has a compressive wavenumber modulation in the wavenumber $g(x)$:

\begin{equation}
u(x,0)=-2g(x)\left(\frac{1-\sqrt 2\cos(sg(x))}{3-2\sqrt 2\cos(xg(x))}\right)
\end{equation}

with 

$$g(x)=\alpha\left(1-\frac{2}{\pi}\arctan \beta x\right)+b_0.$$
\end{remark}

\vspace{0.3cm}
\subsection{The ILW equation}

As recalled in the previous section, the formal framework of inverse scattering for the ILW equation was given in \cite{KSA}, \cite{KAS}. We are not aware of rigorous results using IST to solve the Cauchy problem even for small initial data.

\section{Related results and conjectures}
\subsection{The modified cubic BO and ILW equations}
The modified cubic BO and ILW equations write respectively \footnote{Note that there are different from the modified BO and ILW equations which are linked to the usual BO and ILW equations via a Miura transform, see Subsection 6.6.} 

\begin{equation}\label{mBO}
u_t+u^2u_x-Hu_{xx}=0
\end{equation}

and 

\begin{equation}\label{mILW}
u_t+u^2u_x+\frac{1}{\delta}u_x+\mathcal T(u_{xx})=0.
\end{equation}

Contrary to the modified KdV equation they do not seem to be completely integrable. We will not describe in details the local well-posedness results since most of the used methods are close to that of the original one. We refer for instance to \cite{MR2, KeTa} and the references therein for the local Cauchy problem in $H^s(\R), s\geq \frac{1}{2}$  and to \cite{MR} for the local Cauchy problem in $H^1(\T).$ Scattering results of small solutions can be found in \cite{HN1,HN2}.

Using the factorization technique of Hayashi and Naumkin, Naumkin and Sanchez-Suarez \cite{NaSa} study the large time asymptotics of small solutions of the modified ILW equation equation in suitable weighted Sobolev spaces. Difficulties arise from the non-homogeneity of the symbol. 

 More precisely, the main result in \cite{NaSa} states that for any initial data 

$u_0\in H^3(\R)\cap H^{2,1}(\R)$ with $||u_0||_{H^3\cap H^{2,1}}\leq \epsilon, $ and $\int_\R u_0(s)dx =0,$ there exists a global solution $u\in C([0,\infty); H^3(\R))$ of the modified ILW satisfying

$$|u(t)|_\infty \leq C\epsilon t^{-1/2}.$$
There exist moreover unique modified final states describing the asymptotics of $u(t)$ for large times.

\vspace{0.3cm}
The equations \eqref{mBO} and \eqref{mILW}   are specially  interesting from a mathematical point of view since they are  $L^2$ critical, as is the quintic generalized KdV equation and one thus expects a finite type blow-up phenomena. This has been established for the modified BO equation in \cite{MaPi} as follows.

 Let $Q\in H^{1/2}(\R), Q>0$ be the unique ground state solution of the "elliptic" equation

$$D^1Q+Q=Q^3$$

constructed by variational arguments (see \cite{We, ABS}) and whose uniqueness was established in \cite{FL}. Then it is proven in \cite{MaPi} that there exists a solution $u$ of \eqref{mBO} satisfying $|u(t)|_2=|Q|_2$ and

$$u(t)-\frac{1}{\lambda^{\frac{1}{2}t}}\left(\frac{.-x(t)}{\lambda(t)}\right)\to 0\quad\text{in}\;H^{1/2}(\R)\;\text{as}\;t\to 0$$

where 

$$\lambda(t)\sim t,\; x(t)\sim -|\ln t|\quad \text{and}\quad ||u(t)||_{\dot{H}^{1/2}}\sim t^{-\frac{1}{2}}||Q||_{\dot{H}^\frac{1}{2}}\;\text{as}\; t\;\to 0.$$

This result obviously implies the orbital instability of the ground state. Also this blow-up behavior is unstable since any solution merging from an initial data with $|u_0|_2<|Q|_2$ is global and bounded in $L^2.$ The proof is inspired by a corresponding one for the $L^2$ critical generalized KdV (see \cite{MaMeRa}).There are however  new difficulties  to overcome. First, the slow spatial decay of the solitary wave $Q$ creates serious difficulties to construct a blow-up profile. Next difficulties arise when considering localized versions of basic quantities such as mass or energy. Contrary to the case of the $L^2$ critical KdV equation, standard commutator estimates are not enough and one has to use suitable localization arguments as in \cite{KeMa, KMR}.

\begin{remark}
A similar result is not known but expected to be true for the modified ILW equation \eqref{mILW}.
\end{remark}

\subsection{The periodic case}
In order to maintain a reasonable size, we focussed in the present paper on the Cauchy problem in the real line. For the sake of completeness, we however briefly review here the state of the art for the Cauchy problem posed in the circle $\T$  which leads to many interesting issues, both on the IST and PDE sides.

The periodic BO equation writes

\begin{equation}\label{BOper}
u_t+uu_x-\mathcal Hu_{xx}=0,
\end{equation}

where 

$$\mathcal Hf(x)=\frac{1}{L}PV\int_{-L/2}^{L/2}\cot\left(\frac{\pi(x-y)}{L}\right)f(y)dy,$$

while the periodic ILW equation writes

\begin{equation}\label{ILWper}
u_t+2uu_x+\delta^{-1}u_x-(\mathcal T_\delta u)_{xx}=0, \; \delta >0,
\end{equation}

where

$$\mathcal T_\delta u(x)=\frac{1}{L}PV \int_{-L/2}^{L/2}\Gamma_{\delta,L}(x-y)u(y)dy$$

with

$$\Gamma_{\delta,L}(\xi)=-i\sum_{n\neq 0}\coth\left(\frac{2\pi n\delta}{L}\right)e^{2in\pi\xi/L}.$$

We are not aware of  rigorous results via inverse scattering methods for both the periodic BO and ILW equations, see however the recent progress in \cite{GeKa} for the BO case.

We recall that the formalism of IST for the periodic ILW can be found in \cite{AFSS}. 

\vspace{0.3cm}
On the other hand many results have been  obtained by PDE methods.

The results in \cite{ABFS} still remain true since they do not rely on dispersive properties of the linear group, yielding global well-posedness of the Cauchy problem for the periodic BO or ILW equations in $H^s(\T), s>\frac{3}{2}.$

 Concerning lower regularity results, L. Molinet (\cite{Mol2}) proved the global well-posedness in the energy space $H^{1/2}(\T)$ by combining Tao's normal form and estimates in Bourgain spaces. He also proved that the flow map $u_0\to u(.,t)$ is not uniformly continuous in $H^s(\T), s>0.$ The proof of this fact  is inspired by a corresponding result in \cite{KT} for the Cauchy problem in $H^s(\R), s>0$ but with a simpler proof, using the property that if $u$ is a solution of the BO equation with initial data $u_0$, so is $u(.+\omega t,t)+\omega, \forall \omega \in \R$ with initial data $u_0+\omega.$  Actually this lack of uniform continuity has nothing to do with dispersion and is linked to the Burgers equation (nonlinear transport equation) 
 
 $$u_t+uu_x=0.$$
 Actually the locally Lipschitz property of the flow map is recovered in the space $H^s_0(\T)$ of $H^s(\T)$ functions with zero mean or more generally on hyperplanes of functions with fixed mean value.
 
 This result is improved in \cite{Mol} where  the global well-posedness  in $L^2(\T)$ is proven, the flow map being Lipschitz (in fact real analytic) on every bounded set of $H^s_0(\T), s\geq 0$. Note that uniqueness then holds in a space smaller than $C([0,T]; L^2(\T))$ but containing the limit in $C([0,T]; L^2(\T))$ of smooth solutions of BO. The result is sharp in the sense that the flow map (if it can be defined and coincides with the standard flow-map on $H^\infty_0(\T))$ cannot be $C^{1+\alpha}, \alpha >0,$ from $H^s_0(\T)$ into $H^s_0(\T)$ as soon as $s<0.$
 
 Note that one also finds in \cite{MP2} a simplified proof of the global well-posedness in $L^2(\T)$  with unconditional uniqueness in $H^{1/2}(\T).$
 
 The ill-posedness in $H^s(\T)$ has been reinforced in the main result in \cite{Mol4} which insures that  for all non constant function $\phi\in L^2(\R)$ and all $T>0$, there exists $t\in (0,T)$ such that the map which associates $u_0$ to the corresponding solution of BO $u(t)$  is discontinuous at $\phi$ in any Sobolev space $H^s(\T), s<0.$ This strong ill-posedness property completes exactly the well-posedness result of \cite{Mol}.  This discontinuity result holds also on hyperplanes of functions with a given mean value.
 
\vspace{0.3cm}
As aforementioned, the results in \cite{MV} apply to the ILW equation, yielding global well-posedness in the energy space $H^{1/2}(\T).$ We are not aware  of well-posedness results in $L^2(\T)$   for the ILW equation.
 
 \vspace{0.3cm}
 As in the non-periodic case, both the BO and ILW equations possess  periodic solitary waves. Explicit formulas can be found in \cite {JE, CL, NM} but, as noticed in \cite{AFSS} they apper to be incorrect. An explicit formulation is given in \cite{Mil}.
 The L-periodic traveling wave of the BO equation is in the formulation of  \cite {AN}:
 
 \begin{equation}\label{SWBO}
 \phi_c(x)=\frac{4\pi}{L}\frac{\sinh (\gamma)}{\cosh(\gamma)-\cos(2\pi x/L)},
 \end{equation}
 
 where $\gamma>0$ satisfies $\tanh(\gamma)=\frac{2\pi}{cL},$
  implying that $c>\frac{2\pi}{L}.$
  The profile of an even, zero mean, periodic traveling wave of the ILW equation is given by (see \cite {ACN} and also \cite{NM}) :

  \begin{equation}\label{SWILW}
  \phi_c(x)=\frac{2K(k)i}{L}\left( Z(\frac{2K(k)i}{L}(x-i\delta);k)-Z(\frac{2K(k)i}{L}(x+i\delta);k)\right),
  \end{equation}
  
  where $K(k)$ denotes the complete elliptic integral of first kind, $Z$ is the Jacobi zeta function and $k\in (0,1).$ Note that for fixed $L$ and $\delta,$  the wave speed c and the elliptic modulus k must satisfy specific constraints.

We also refer to \cite{Par} for another derivation and a clear and complete discussion of the periodic solutions of the ILW equation.
  
\vspace{0.3cm}
The orbital stability of the periodic traveling waves is proven in \cite{AN} for the BO equation and in \cite{ACN} for the ILW equation.

\vspace{0.3cm}
A program (see \cite{Tz, Deng, TV1, TV2, TV3, DTV}) has been devoted to the construction of an infinite sequence of invariant measures of Gaussian type associated  to the conservation laws of the Benjamin-Ono equation.

 Those invariant measures $\lbrace \mu_n\rbrace$ on $(L^2,\phi_t)$ where $\phi_t$ is the BO flow on $L^2(\T)$ satisfy

$$ \mu_n \;\text{is concentrated on}\; H^s(\T)\;\text{for}\; s<n-\frac{1}{2},$$

$$\mu_n(H^{n-1/2}(\T)) =0.$$

$\mu_n$ is formally defined as a renormalization of 

$$d\mu_n(u)=e^{-E_n(u)}du=e^{R_n(u)}e^{||u||_n^2}du,
$$
where  we have denoted 

$$E_n(u)=||u_n||^2+R_n(u)$$ 

the sequence of conservation laws of the BO equation, $||.||_n$ denoting the homogeneous Sobolev norm of order $n$ and $R_n$ a lower order term.

These results do not apply to infinitely smooth solutions since $\mu_n(C^\infty(\T))=0$ for all $n.$ Sy \cite{Sy} has constructed a measurable dynamical system for BO on the space $C^\infty(\T)$, namely  a probabilistic measure $\mu$ invariant by the BO flow, defined on $H^3(\T))$ and satisfying in particular 
$\mu(C^\infty(\T))=1.$
\begin{remark}
No similar results seem to be available for the ILW equation.
\end{remark}

\subsection{Zero dispersion limit}

Interesting issues arise when looking at the zero dispersion limit of the BO equation, that is the behavior of solutions to 

\begin{equation}\label{BOeps}
u_t+uu_x-\epsilon Hu_{xx}=0,\; \epsilon >0
\end{equation}

when $\epsilon \to 0.$

The corresponding problem for the KdV equation has been first extensively studied by Lax and  Levermore  in the eighties \cite{LaLi}, see also \cite{Ven} and the survey article \cite{Miller}. For the BO equation, as for the KdV equation, at the advent of a shock in the dispersionless Burgers equation, the solution of \eqref{BOeps} is regularized by a dispersive shock wave (DSW). In the DSW region, the solution of \eqref{BOeps} may be formally approximated using Whitham modulation theory, see \cite{Wh}. Unlike the case of the KdV equation the modulation equations for the BO equation are fully uncoupled (see \cite{DK}), consisting of several independent copies of the inviscid Burgers equation.

A formalism for matching the Whitham modulation approximation for the DSW onto inviscid Burgers equation in the domain exterior to the DSW was developed by Matsuno \cite{Mat14, Mat15} and Jorge, Minzoni and Smyth in \cite{JMS}.

Recent rigorous results have been obtained by P. Miller and co-authors. Miller and Xu \cite{MX} partially confirmed the above formal results by computing rigorously the weak limit (modulo an approximation of the scattering data) of the solution of the Cauchy problem for \eqref{BOeps} for a class of positive initial data, using IST techniques, and by developing  an analogue for the BO equation of the method of Lax-Levermore for the KdV equation.

Using the exact formulas for the scattering data of the BO equation valid for general rational potential with simple poles that they obtained in \cite{MW2}, Miller and Wetzel \cite{MW} analyzed rigorously the scattering data in the small dispersion limit, deducing in particular precise asymptotic formulas for the reflection coefficient, the location of the eigenvalues and their density, and the asymptotic dependence of the phase constant associated with each eigenvalue on the eigenvalue itself.  Such an analysis seems to be unknown for more general potentials.

We refer to \cite{MX2} for a study  of the BO hierarchy with positive initial data  in a suitable class  in the zero-dispersion limit. We also mention   \cite{ENS} for a description of dispersive shock waves for a class of nonlinear wave equations with a nonlocal BO type dispersion which does not use the integrability of the equation.

We are not aware of similar small dispersion limit results for the ILW equation.

\subsection{The soliton resolution conjecture}

Since IST methods provide so far rigorous results on the Cauchy problem only for small enough initial data, one cannot  use  them to prove the soliton resolution conjecture (see \cite{Tao}) as in the KdV case (see {\it eg} \cite{Schu} and the references therein). However we strongly believe that it holds true for the BO and ILW equations, as suggested for instance by the following numerical simulations in \cite{KS2} (see also \cite{PD, MPST} for other illuminating simulations and \cite{PD2} for a theoretical analysis of the spectral method in \cite{PD}). \footnote{ Note that the soliton resolution conjecture might also be valid  for the Benjamin equation \eqref{Benjamin} which is not integrable.}

\begin{figure}[htb!]
 \includegraphics[width=\textwidth]{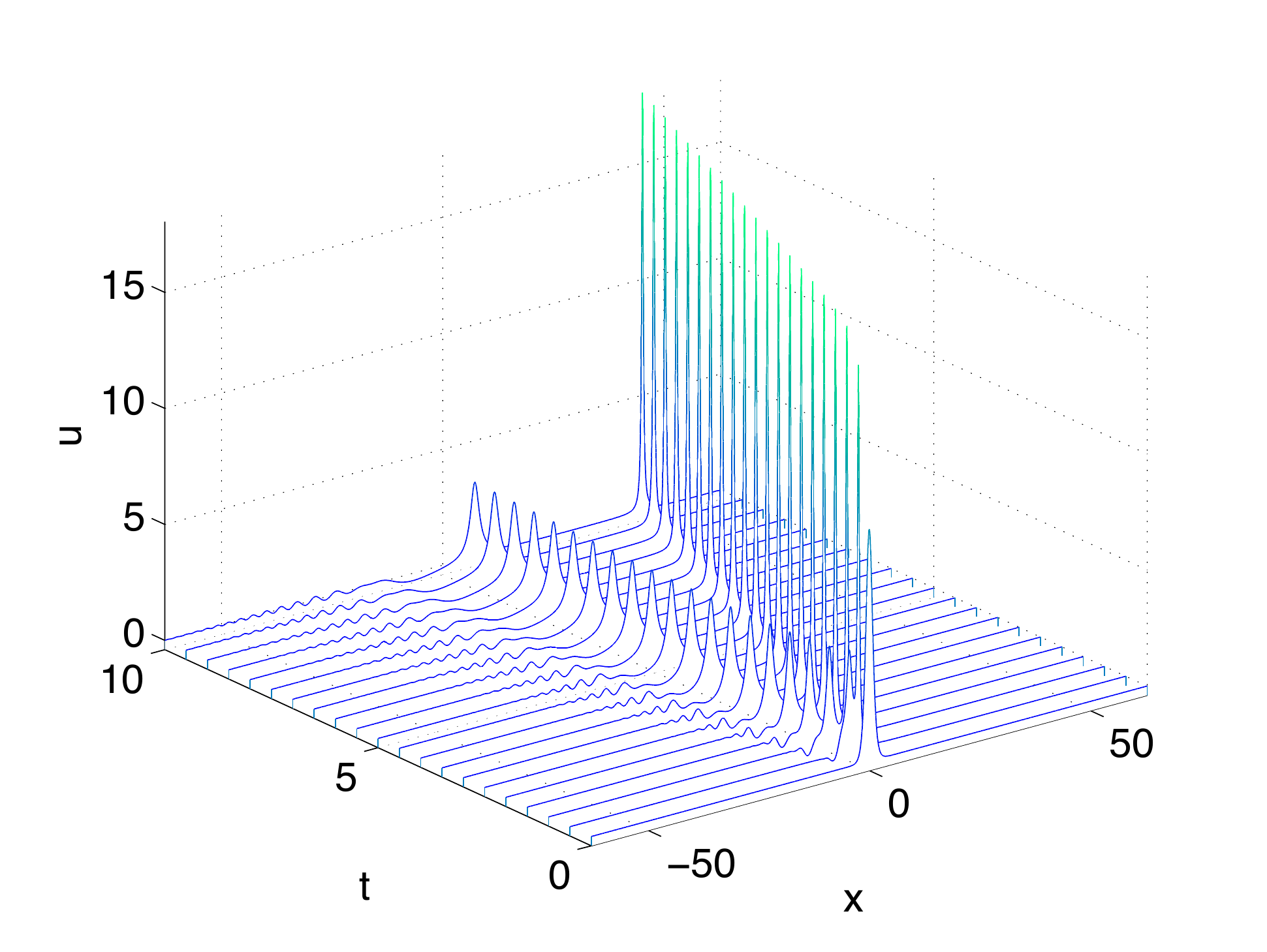}
 \caption{Solution to the BO equation for the initial data $u_0=10 \text{sech}^2 x$} 
 \label{BO10sech}
\end{figure}

We first show the formation of solitons from localized initial data for the BO equation in 
Fig.~\ref{BO10sech}. Note a tail of dispersive oscillations 
propagating to the left. 

\begin{figure}[htb!]
 \includegraphics[width=\textwidth]{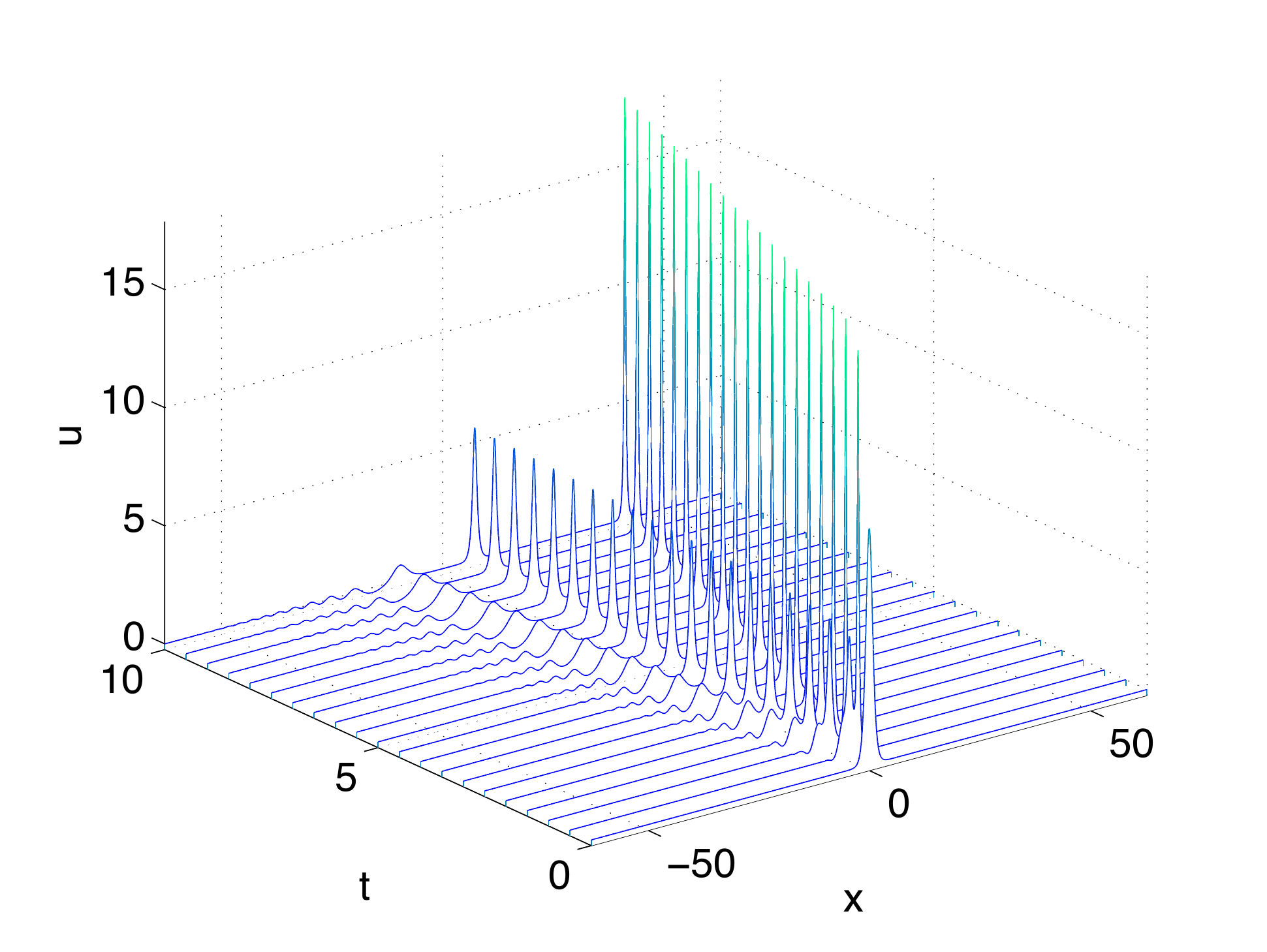}
 \caption{Solution to the ILW equation with $\delta=1$  for the initial data $u_0=10 \text{sech}^2 x$}
 \label{IWL10sech}
\end{figure}

A similar  decomposition of localized initial data into solitons and 
radiation is shown for the ILW equation in Fig.~\ref{IWL10sech}. 
Note that this case is numerically 
easier to treat with Fourier methods since the soliton solutions are 
more rapidly decreasing (exponentially instead of algebraically) than for the fKdV, fBBM and BO equations. The different shape of the solitons is also noticeable in 
comparison to Fig.~\ref{BO10sech}. 

More generally, we conjecture that some form of soliton resolution holds for the class of fractional KdV equations

\begin{equation}\label{fKdV}
u_t+uu_x-|D|^\alpha u_x=0,
\end{equation}

when $\alpha >\frac{1}{2},$ the Cauchy problem being then globally well-posed in the energy space $H^{\alpha/2}(\R).$ 
\footnote{Recall that the ground state solution (that exists and is unique when $\alpha>\frac{1}{3}$) is orbitally stable if and only if $\alpha>\frac{1}{2}$, see \cite{LPS3}.}  first step towards this conjecture is achieved in \cite{MPV} where global well-posedness is established in the energy space when $\alpha >\frac{6}{7}.$

The  simulation in Figure 3 below for the fKdV with $\alpha =0.6$ suggests the decomposition into solitons plus radiation.
\begin{figure}[htb!]
 \includegraphics[width=\textwidth]{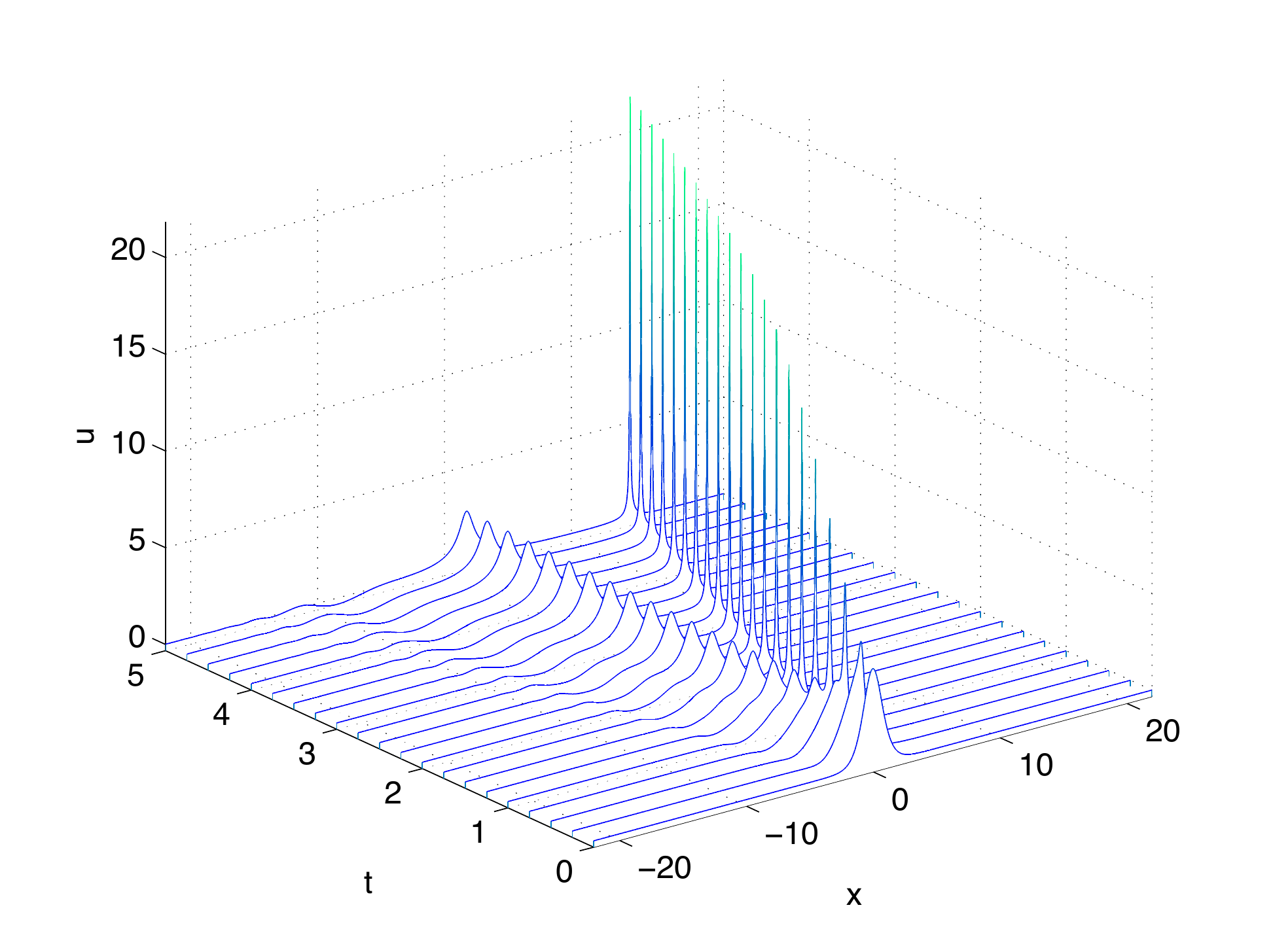}
 \caption{Solution to the fKdV equation with $\alpha=0.6$ for the initial data $u_0=5\text{sech}^2x$ }
 \label{gBO5sechalpha06water}
\end{figure}

\section{Varia}

We describe here a few more qualitative results on the Benjamin-Ono and the Intermediate Long Wave equations and some of their natural extensions (higher order equations, two-dimensional versions...).

\vspace{0.3cm}
\subsection{Damping of solitary waves}
In the real physical world, various dissipative mechanisms affect the propagation of dispersive waves. The dissipative term can be local {\it eg} 
$-\partial_{xx}$ or non local, see for instance \cite{OS, OS2, KakMa, ER} for physical examples. The Cauchy problem for the dissipative BO or ILW equations is studied in \cite{JCS} and more recently in \cite{MV, Mol3}. Of course the dissipative term destroys the Hamiltonian structure and the complete integrability. The solutions of the Cauchy problem then tends to zero as $t\to \infty$ (see \cite{BoLu, Dix, Dix2} and the references therein for a description of the large time asymptotics of solutions of the Benjamin-Ono-Burgers and related  equations).

However we are not aware of rigorous results on  the quantitative influence of dissipation in the solitonic structure of the BO, or ILW equation \footnote{Actually this question seems to be open also for the KdV equation}.

In \cite{MSKK} the Whitham method is applied to the Benjamin-Ono-Burgers equation 

\begin{equation}\label{BOB}
u_t+uu_x-Hu_{xx}=\epsilon u_{xx}, \epsilon \ll 1
\end{equation}
to give a formal description of the solution of the Cauchy problem corresponding to a steplike initial data as $\epsilon \to 0$ .

\vspace{0.3cm}
\subsection{Weighted spaces} The Cauchy problem for the BO equation in weighted spaces is studied in \cite {Io2, HKO, Flo, FoFu, FoLi, FLP1,FLP2} yielding smoothing and unique continuation properties \cite{Io3}. A result on the propagation of regularity for  solutions to the Cauchy problem  associated to the Benjamin-Ono equation is proved in \cite{ILP}. Essentially, if $u_0\in H^{3/2}(\R)$ is such that its restriction belongs to $H^m(b,\infty)$ for some $m\in \mathbb{Z}, m\geq 2$ and some $b\in \R$, then  the restriction of the corresponding solution  $u(.,t)$ belongs to $H^m(\beta, \infty)$ for any $\beta \in \R$ and any $t>0.$ This shows that the regularity of the datum travels to the left with infinite speed.

Another intesting result is in \cite{ILS} where the Cauchy problem for BO is solved with bore-like initial data.

\vspace{0.3cm}
\subsection{Control}  Control and stabilization issues were investigated by Linares and Rosier, \cite{LinRos}, in $H^s(\T), s>\frac{1}{2}$, and in \cite{LLR} in $L^2(\T)$ for the periodic BO equation.

A typical exact controllability  (via a forcing term) result from \cite{LLR} is as follows. In order to keep the mass conserved, the control input is chosen in the form 

$$ (\mathcal G h)(x,t)=a(x)\left (h(x,t)-\int_\mathbb{T} a(y)h(x,t)dy\right ),$$

where $a$ is a given smooth nonnegative function such that $\lbrace x\in \mathbb{T}; a(x)>0\rbrace =\omega$ and with mass one.

\begin{theorem}
(i) (Small data) For any $T>0$ there exists some $\delta >0$ such that for any $u_0, u_1\in L^2(\mathbb{T})$ with 

$$|u_0|_2\leq \delta,\; |u_1|_2\leq \delta\;\text{and}\; \int_{\mathbb{T}} u_0=\int_{\mathbb{T}} u_1,$$

one can find  a control input $h\in L^2(0,T;L^2(\mathbb{T})$ such that the solution $u$ of the system

\begin{equation}
    \label{abcd}
    \left\lbrace
    \begin{array}{l}
    u_t-Hu_{xx}+uu_x=\mathcal G h \\
    u(x,0)=u_0(x),
\end{array}\right.
    \end{equation}
satisfies $u(x,T)=u_1(x)$ on $\mathbb{T}.$

(ii) (Large data) For any $R>0$ there exists a positive $T=T(R)$ such that the above property holds for any$u_0, u_1\in L^2(\mathbb{T})$ with

$$|u_0|_2\leq R,\; |u_1|_2\leq R\; \text{and}\; \int_{\mathbb{T}} u_0=\int_{\mathbb{T}} u_1$$
\end{theorem}

\begin{remark}
1. We refer to \cite{LLR} for stabilizations results.

2. Related results for the {\it linear} problem were obtained in \cite{LO}.
\end{remark}

This result relies strongly on the bilinear estimates proved in \cite{MP2}.

  We do not know of corresponding results for the ILW equation.

\vspace{0.3cm}
\subsection{Initial-boundary value problems}

The global well-posedness of the initial- boundary value problem the BO and the ILW equation on the half-line with zero boundary condition at $x=0$ is proven respectively in \cite{HK} and in \cite{AK}, for small initial data in suitable weighted Sobolev space. Moreover the  long time asymptotic is given, for instance, one gets for the ILW equation \cite{AK}  

$$u(x,t)=\frac{1}{3\pi \delta t}\frac{x}{(\sigma t)^{1/3}}\text {Ai}\left( \frac{x}{(\delta t)^{1/3}}\right)
+\min\left( 1, \frac{x}{\sqrt t}\right)O\left(t^{-1-a/2}\right),$$
where $a>0$ and $\text{Ai}$ is the Airy function.

\vspace{0.3cm}
\subsection{Transverse stability issues}

The Kadomtsev-Petviashvili (KP) equation was introduced heuristically in \cite{KaPe} to study the 
transverse stability of the KdV soliton with respect to long, weakly transverse perturbations. Since their formal analysis does not depend on the dispersive term in the KdV equation, it applies as well as to the BO and ILW equations, yielding the KP-II version of those one dimensional equations, that is, respectively

\begin{equation}\label{KP-BO}
u_t+u_x+uu_x-Hu_{xx}+\partial_x^{-1}u_{yy}=0,
\end{equation}

and

\begin{equation}\label{KP-ILW}
u_t+uu_x+\frac{1}{\delta}u_x+\mathcal T(u_{xx})+\partial_x^{-1}u_{yy}=0
\end{equation}

Those equations (that are not known to be integrable) were in fact derived formally in the context of internal waves  in \cite{AbSe, ChCa, GZ, Mat12}. 

Note however that they suffer from two shortcomings of the usual KP-II equation. First, as was, noticed in \cite{Lannes1}, the (artificial) singularity at $\xi_1=0$ of the symbol $\xi_1^{-1}\xi_2^2$ of $\partial_x^{-1}\partial_{yy}$ induces a poor error estimate when comparing the KP-II equation with the full water wave system in the appropriate regime. Roughly speaking, one get (see \cite{Lannes1} and \cite{LS}) an error estimate of the form, in suitable Sobolev norms:

$$||u_{KP}-u_{WW}||=o(1),$$

the "correct" error should be, as in the KdV case,  $O(\epsilon^2 t)$ where $\epsilon$ is the small parameter measuring the weak nonlinear and long wave effects. It is very likely that such a poor precision also holds for the KP-BO and KP-ILW equations.

The second shortcoming concerns the zero mass in $x$ constraint also arising from the singularity of the symbol (see \cite{MST1}). This of course do not exclude that those KP type equations predict at least qualitatively  the behavior of real waves, see {\it eg} \cite{AB,Se}.

 Coming back to rigorous results, the local Cauchy problem (for a larger class of equations including KP-BO and KP-ILW) is studied in \cite{LPS2} and the global existence and scattering of small solutions is proven in \cite{H-GM}. 

It was proven in \cite{LPS3} that \eqref{KP-BO} does not possess any non trivial (lump type) {\it localized} solitary wave. The proof extends as well to \eqref{KP-ILW}.

Recall that the KdV soliton is transversally $L^2$ stable with respect to weak transverse perturbations described by the KP II equation (\cite{MT, Mi}). This result is unconditional since  it was established in  \cite{MST5} that the Cauchy problem for KP II is globally well-posed in $H^s(\R\times \T), s\geq 0,$ or for all initial  data of the form $u_0+\psi_c$ where $u_0\in H^s(\R^2), s\geq 0$ and $\psi_c(x-ct,y)$ is a solution of the KP-II equation such that for every $\sigma \geq 0,$ $(1-\partial_x^2-\partial_y^2)^{\sigma/2}\psi_c$ is bounded and belongs to $L_x^2L_y^\infty(\R^2).$ Note that such a condition is satisfied by the value at $t=0$ of any soliton or N-soliton of the KdV equation. 

A priori this condition is not satisfied by a function $\psi$ that is not decaying along a line $\lbrace (x,y); x-vy=x_0\rbrace$  such as for instance a KP-II line soliton that writes $\psi(x-vy-ct).$ However, one checks that the change of variable 
\begin{equation}\label{change}
X=x+vy-v^2t, Y=y-vt, T=t
\end{equation} 
leaves the KP-II equation invariant and thus the global well-posedness results hold for an initial data  that is  localized perturbation of the KP-II line soliton. Observe that the same transformation leaves also invariant the KP-BO and KP-ILW equations.

It is natural to conjecture that a similar stability result holds true for the  BO and ILW soliton with respect to \eqref{KP-BO} and \eqref{KP-ILW}. We refer to \cite{AbSe} for formal arguments in favor of this conjecture.

Such a result would be however a conditional one since no global existence result is available for both  \eqref{KP-BO} and \eqref{KP-ILW}, 

It is worth noticing however that it was observed in \cite{SM}, following the method of \cite{KSF} that the {\it periodic} BO solitary waves could be spectrally transversally unstable for some ranges of frequencies. This deserves further investigations. 

\vspace{0.3cm}
\subsection{Modified BO and ILW equations}

It is well known that a real solution of the focusing modified KdV equation

$$v_t+6v^2v_x+v_{xxx}=0$$

is sent by the Miura transform

$$u=v^2+iv_x$$

to a complex solution of the KdV equation 

$$u_t+6uu_x+u_{xxx}=0.$$

It turns out that a similar fact holds for the BO and ILW equations, see \cite{Nak3, Nak4, STA, Mat3} and specially \cite{SC} to which we refer for an explanation of the mathematical origin of the MILW equation. 

The modified ILW and BO equations  write respectively

\begin{equation}\label{MILW}
v_t+\beta v_x(e^v-1)+\frac{1}{\delta}v_x+v_x\mathcal T(v_x)+\mathcal T(v_{xx})=0,
\end{equation}

\begin{equation}\label{MBO}
q_t+\alpha q_x(e^q-1)+q_xH(q_x)+H(q_{xx})=0,
\end{equation}

where $\alpha, \beta$ are real constants.

\begin{remark}
We are not aware of results on the Cauchy problem for \eqref{MBO} or  \eqref{MILW} by PDE methods.
\end{remark}

Both  equations  \eqref{MBO} and \eqref{MILW} are integrable soliton equations. In particular the MILW possesses an infinite number of conserved quantities, \cite{STA}, a linear scattering theory \cite{STA}, a B\"{a}cklund transform \cite{STA, Nak3} and multi-soliton solutions, \cite{Nak3}. The formal IST for the MBO and MILW is studied in \cite {SC2, SC} respectively.

One can check, see \cite{SC} that the Miura type transformation

$$u=\frac{1}{2}\lbrace \mathcal T (v_x)+\beta(e^v-1)+iv_x\rbrace$$

maps real-valued solutions of \eqref{MILW} with $\alpha=-\beta$ into complex-valued solution of the ILW equation.

A similar transformation holds in the BO case, see \cite{SC2}.

\vspace{0.3cm}
\subsection{Higher order  BO and ILW}
Higher order BO and ILW equations (as higher order KdV) appear in two different contexts.

 1. First when going to next orders  in the expansion of the two-layer system in the ILW /BO regime, see Introduction. For instance,  one gets in the notations of \cite{CGK} (see also \cite{Mat12} for a similar equation)   the next equation in the BO regime

\begin{multline}\label{BO4}
u_t=\frac{\rho h_1^2}{2\rho_1^2}A^2 H(u_{xx})+\frac{3\sqrt 2}{4\rho_1}Auu_x\\-\frac{\sqrt 2}{2}\epsilon\frac{\rho h_1^2}{\rho_1^2}
A[\partial_x(uH(u_x)+H\partial_x(uu_x)] +\frac{\epsilon}{2}\left(\frac{\rho^2h_1^2}{\rho_1^2}-\frac{h_1^2}{3}\right)A^2u_{xxx},
\end{multline}

where   $A= \left(\frac{g\rho_1(\rho-\rho_1)}{h_1}\right)^{1/2},$ and where $\rho>\rho_1$ are the densities of the two fluid layers, $h_1$ is the height of the upper layer, $g$ the constant of gravity and $\epsilon >0$ is a small parameter.  

This process could actually be continued to obtain an infinite sequence of higher order BO equations. It applies as well to the ILW equation, yielding the following equation, in the above notations of \cite{CGK}, where here $\mathcal T_h$ is the Fourier multiplier with symbol $-i\coth(\epsilon h|\xi|),$ $h$ being the depth of the lower layer:

\begin{multline}\label{ILW4}
u_t=\frac{\rho h_1^2}{2\rho_1^2}A^2 \mathcal T_h(u_{xx})+\frac{3\sqrt 2}{4\rho_1}Auu_x-\\\frac{\sqrt 2}{2}\epsilon\frac{\rho h_1^2}{\rho_1^2}
A[\partial_x(u\mathcal T_h(u_x)+\mathcal T_h\partial_x(uu_x)] -\frac{\epsilon}{2}\left(\frac{\rho^2h_1^2}{\rho_1^2}(\mathcal T_h)^2+\frac{h_1^2}{3}\right)A^2u_{xxx}.
\end{multline}

An equation like \eqref{BO4} was also derived in \cite{Mat18} where a solitary-wave solution of the equation is obtained by means of a singular perturbation method. The characteristics of the solution are discussed in comparison with those for a higher-order BO equation of the Lax type, see next section.

\vspace{0.3cm}
2. Second, one gets a hierarchy of (integrable)  higher order BO or ILW equations by considering the Hamiltonian flows of the successive conservation laws, see {\it eg} \cite{Mat3}.

For instance, the next equations in the hierarchy are for the ILW and BO equation   given respectively by

$$u_t=\partial_x\text{grad}\;I_4(u)$$

and

$$u_t=\partial_x\text{grad}\;J_4(u),$$

\vspace{0.3cm}
where $I_4$ and $J_4$ were defined in Chapter 2, that is 

\begin{equation}\label{HOBO}
u_t-3u^2u_x+4u_{xxx} +3(uHu_x)_x-3H(uu_x)_x=0
\end{equation}

and

  \begin{equation}\label{HOILW}
u_t-3u^2u_x-3(u\mathcal T u_x)_x-3\mathcal T(u_{xx})+u_{xxx}-3\mathcal T^2(u_{xxx})-\frac{1}{\delta}[18uu_x+9\mathcal T u_{xx}]-\frac{3}{\delta^2}u_x=0
\end{equation}

We are not aware of rigorous results on \eqref{HOBO}, \eqref{HOILW} by IST methods. 

On the other hand, it was shown in \cite{Pi} that \eqref{BO4} is quasilinear in the sense that, as for the BO equation,  the flow map cannot be $C^2$ in any Sobolev space $H^s(\R), s\in \R$  (the proof applies as well to \eqref{HOBO}).  In \cite{LPP} the Cauchy problem for \eqref{BO4} is proven to be locally well-posed in $H^s(\R), s\geq 2$  and in some weighted Sobolev spaces.

This result was significantly improved in \cite{MP} where the global well-posedness of \eqref{BO4} was established in $H^s(\R), s\geq 1,$ and thus in particular in the energy space $H^1(\R)$. The main difficulty is of course to obtain the local well-posedness in $H^1(\R).$ This is achieved by introducing, as in \cite{Tao}, a gauge transformation that weakens the high-low frequency interactions in the nonlinear terms. Such a transformation was already used in \cite{LPP} to obtain the $H^2$ well-posedness but it is combined here with a Besov version of Bourgain spaces, together with the full Kato smoothing effect for functions that are localized in space frequencies.

Another result in \cite{MP} concerns the limiting behavior of solutions $u_\epsilon$ when $\epsilon \to 0.$ A direct standard compactness method (as used for instance in the BO-Burgers equation in the zero dissipation limit) does not seems to work since the two leading terms (of order one) in the Hamiltonian  $I_4$ have opposite signs. However it is shown that the solution of \eqref{BO4} converges in $L^\infty(0,T;H^1(\R)), \forall T>0$ to the corresponding solution of BO as $\epsilon \to 0$  provided the ratio of density $\frac{\rho_1}{\rho}$ equals $\frac{1}{\sqrt 3}.$

This shortcoming of \eqref{BO} might not occur when using the method of \cite{BLS} to derive asymptotic models of internal waves since it is more flexible and leads to (equivalent in the sense of consistency) {\it families} of models. In particular at least one member of the family of higher order BO equations might contain one for which the limit to BO holds, for any ratio of densities.

As mentioned in \cite{MP}, it is likely that the previous results hold true for the higher order ILW equation \eqref{ILW4}.

We are not aware of similar results for the next equations in the BO or ILW hierarchy \eqref{HOBO} and \eqref{HOILW} but it is likely that the method in \cite{MP} could be applied. 

For any of those higher order equations (\eqref{HOBO}, \eqref{HOILW}, \eqref{BO4}, \eqref{ILW4}) questions about existence and stability of solitary wave solutions seem to be widely open except for the Lax hierarchy of the BO equation for which Matsuno \cite{Mat7} has established the Lyapunov stability of the N-soliton  in the hierarchy, using in particular results from the IST. One should also mention \cite{Mat18} where by using a multisoliton perturbation theory it is proven analytically that the overtaking collision between two-solitary waves exhibits a phase shift, the amplitudes being not altered after interaction.

\vspace{0.3cm}
\subsection{Interaction of solitary waves}

As noticed in \cite{WP} where one can find a nice overview of the subject and a  fascinating analysis of observed (via imaging technique)  oblique wave-wave interactions on internal waves in the strait of Georgia, "although nonlinear interactions that occur when two large internal waves collide at oblique angles are often observed in the natural world, quantitative and theoretical aspects  of these interactions are only poorly understood".

This  important topic  has been first extensively studied for shallow water gravity waves, mainly in the context of the KP-II equation.
Although KP-II has no fully localized solitary wave solutions \cite{deBS}, it has a large variety of two-dimensional exact solutions, those wave patterns being generated by nonlinear interactions among several obliquely propagating solitary waves. In particular the resonant interactions among those  solitary waves play a fundamental role in multidimensional wave phenomenon, leading to very complicated patterns. We refer for instance to \cite{BMOT, Se, ChaSat} and to the book \cite{Kodama} and the references therein for an extensive description of those waves. We also refer to \cite{D-Wu} for a study of the direct scattering theory  for KP II on the background of a line soliton. 

It is worth noticing that some of those  patterns look very much like real observed waves (\cite{AB}) and the classical picture below of interaction of line solitons on the Oregon coast.



 
\begin{figure}[htbp]\includegraphics[]{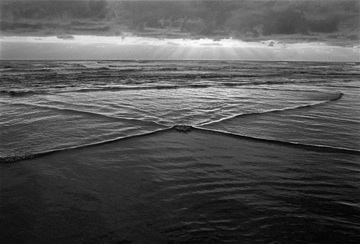}
\end{figure}



According to \cite{Miles, Miles2} the oblique interaction of solitary waves propagating at different directions ${\bf n}_1$ and ${\bf n}_2$ can be classified into two types. The first one occurs when the two solitary waves propagate in almost the same direction ( $1-{\bf n}_1.{\bf n}_2\gg \epsilon$ where $\epsilon$ is the ratio between a typical amplitude of the waves and a typical depth of the fluid) ) and interact for a relatively long time (strong interaction). The second type corresponds to the interaction of solitary waves propagating in almost opposite directions ($1-{\bf n}_1.{\bf n}_2\simeq O(\epsilon))$ and hence the interaction is relatively short (weak interaction). 

Typically the two waves evolves according to their own one-dimensional equation (KdV for surface waves ILW or BO for internal waves) in the case of weak interactions and according to a KP-II like equation in the case of strong interactions.


Such exact multiline solitons are not known to exist for the BO and ILW versions of the KP equation since they are not integrable. Nevertheless it makes prefect sense to look for oblique interactions of internal line solitons. Most of the existing work consists in numerical simulations.


For internal waves the question of oblique interaction of solitary waves has been addressed in particular by Oikawa and Tsuji,  in the context of BO (infinite depth) \cite{OT, TO} and ILW (finite depth \cite{ TO2}). One finds in \cite{TO} numerical simulations of the strong interaction of nonlinear  long waves whose propagation directions are very close to each other. Two initial settings are considered, first a superposition of two BO solitons with the same amplitude and with different directions, and the second one is an oblique reflection of a BO soliton at a vertical wall. It is observed that the Mach reflection does occur for small incident angles and for some incident angles very large stem waves (see below) are observed.

The case of finite depth is considered in \cite{TO2}  under the assumption of a small but finite amplitude. When the angle $\theta$  between the wave normals of two solitons is not small, it is shown by a perturbation method that in the lowest order of approximation the solution is a superposition of two ILW solitons and in the next order of approximation the effect of the interaction appears as position phase shifts and as an increase in amplitude at the interaction center of two solitons. When $\theta$  is small, it is shown that the interaction is described approximately by \eqref{KP-ILW}. By solving it numerically for a V-shaped initial wave that is an appropriate initial value for the oblique reflection of a soliton due to a rigid wall, it is shown that for a relatively large angle of incidence $\theta$ the reflection is regular, but for a relatively small $\theta$  the reflection is not regular and a new wave called stem is generated. The results are also compared with those of the Kadomtsev-Petviashvili (KP) equation and of the two-dimensional Benjamin-Ono equation \eqref{KP-BO}.

\begin{remark}
Recall that BO or ILW equations have the limitation of weak nonlinearity. Other internal waves models without the assumption of weak non-linearity have beed derived in \cite {ChCa2, Mat12}. The system derived in \cite{ChCa2} has solitary wave solutions which are in good agreement with the experimental results in \cite{KB}.

On the other hand, we recall that Matsuno in particular, (\cite{Mat12}), has extended the (formal) analysis leading to  the BO equation to next order in amplitude and derived a higher order equation similar to a higher order BO equation, namely, say for the wave elevation:
\begin{multline}\label{HOMatsuno}
\zeta_t+\zeta_x+\frac{3}{2}\alpha \zeta\zeta_x+\frac{1}{2}\alpha \delta H\zeta_{xx}-\frac{3}{8}\alpha^2\zeta^2\zeta_x\\
+\frac{1}{2}\sigma\alpha\delta[\frac{5}{4}\zeta H\zeta_{xx}+\frac{9}{4}H(\zeta \zeta_x)_x]-\frac{3}{8}(\sigma^2-\frac{4}{9})\delta^2\zeta_{xxx}=0,
\end{multline}

where here $\sigma=\frac{\rho_1}{\rho_2}$ is the ratio of densities, $\alpha=\frac{a}{h}$ the ratio of a typical amplitude of the wave over the depth of the upper layer and $\delta=\frac{h}{l}$ is the ratio of $h$ over a typical wavelength.

The same issues of transverse stability of 1D- solitary waves or of interaction of line solitons arise for those higher order models but have not be addressed yet as far as we know.
\end{remark}

\vspace{0.3cm}
The situation is different when  considering the interaction of long internal gravity waves propagating on say two neighborhood pycnoclines in a stratified fluid, see \cite{LPK, LKK} where a system of coupled ILW equations is derived. We refer to \cite{Gear, GG} for the shallow-water regime, leading to a system of coupled KdV or KP-II type equations.

For  deep water systems, Grimshaw and Zhu \cite{GZ} have shown that in the strong interaction case, each wave is governed by its own ILW equation (BO equation in the infinite depth case), the main effect being a phase shift of order $O(\epsilon),$ $\epsilon$ measuring the wave amplitude. In the case of weak interactions,  and when $\Delta_1, \Delta_2$ are of order $O(\epsilon),$ the interaction is governed by a coupled two-dimensional KP-II type ILW or BO equations.  Here $\Delta_1=|c_m/c_n-\cos \delta|, \Delta_2=|c_n/c_m-\cos \delta|,$ where $\delta$ is the angle between the two directions of propagation and  $c_m, c_n$ are the linear, long wave speeds for wave with mode numbers $m, n.$
We refer to \cite{LPS3} for a mathematical study of those systems.

\vspace{0.3cm}
In any case the rigorous mathematical  analysis of oblique interactions of internal waves deserves further investigations.

\vspace{0.3cm}
\begin{merci}
The Author thanks heartfully the co-organizers of the program, Peter Miller,  Peter Perry and Catherine Sulem for their involvement in the success of the Program. He also thanks  Felipe Linares and Didier Pilod for very useful comments on a preliminary draft. \end{merci}

\bibliographystyle{amsplain}

\end{document}